\documentclass[11pt]{amsart}

\usepackage{amsthm, amsfonts, amssymb, amscd, rotating}
\usepackage[pagebackref,colorlinks]{hyperref}
\usepackage{tikz-cd}
\usepackage{geometry}
\usepackage{marginnote}

\theoremstyle{definition}
\newtheorem{ntn}{Notation}[section]

\theoremstyle{plain}
\newtheorem{lem}[ntn]{Lemma}
\newtheorem{prp}[ntn]{Proposition}
\newtheorem{thm}[ntn]{Theorem}
\newtheorem{cor}[ntn]{Corollary}

\theoremstyle{definition}
\newtheorem{question}[ntn]{Question}
\newtheorem{rem}[ntn]{Remark}
\newtheorem{exa}[ntn]{Example}

\numberwithin{equation}{section}

\newcommand{\N}{\mathbb{N}}
\newcommand{\z}{\mathbb{Z}}
\newcommand{\q}{\mathbb{Q}}
\newcommand{\R}{\mathbb{R}}
\newcommand{\C}{\mathbb{C}}

\newcommand{\F}{\mathbb{F}}

\newcommand{\EE}{\mathcal{E}}
\newcommand{\Aa}{\mathcal{A}}
\newcommand{\BB}{\mathcal{B}}
\newcommand{\GG}{\mathcal{G}}

\newcommand{\WW}{\mathcal{W}}
\newcommand{\TT}{\mathcal{T}}

\newcommand{\PP}{\mathcal{P}}

\newcommand{\KK}{\mathcal{K}}

\newcommand{\mmm}{\mathfrak{m}}
\newcommand{\ep}{\varepsilon}

\newcommand{\ff}{{F^\times}}

\renewcommand{\aa}{{A^\times}}

\newcommand{\of}{{\overline{F}}}
\newcommand{\off}{{\overline{F}^\times}}

\newcommand{\tors}{{{\rm Tor}_1^{\z}}}
\newcommand{\exts }{{{\rm Ext}_{\z}^1}}

\newcommand{\half}{{\Big[\frac{1}{2}\Big]}}

\newcommand{\mt}{\mapsto}
\newcommand{\lan}{\langle}
\newcommand{\ran}{\rangle}
\newcommand{\se}{\subseteq}
\newcommand{\arr}{\rightarrow}
\newcommand{\larr}{\longrightarrow}
\newcommand{\harr}{\hookrightarrow}
\newcommand{\two}{\twoheadrightarrow}

\newcommand{\stabe}{{\rm Stab}}
\newcommand{\GL}{{\rm GL}}
\newcommand{\PGL}{{\rm PGL}}
\newcommand{\PSL}{{\rm PSL}}

\newcommand{\PB}{{\rm PB}}
\newcommand{\PT}{{\rm PT}}
\newcommand{\Ee}{{\rm E}}

\newcommand{\SL}{{\rm SL}}

\newcommand{\Ind}{{\rm Ind}}
\renewcommand{\char}{{\rm char}}
\newcommand{\diag}{{\rm diag}}
\newcommand{\coker}{{\rm coker}}
\newcommand{\im}{{\rm im}}
\newcommand{\ind}{{\rm ind}}

\newcommand{\inc}{{\rm inc}}
\newcommand{\id}{{\rm id}}

\newcommand{\GE}{{\rm GE}}
\newcommand{\PGE}{{\rm PGE}}

\newcommand{\ab}{{\rm ab}}
\newcommand{\St}{\mathrm{St}}
\newcommand{\GSt}{\mathrm{GSt}}
\newcommand{\PGSt}{\mathrm{PGSt}}
\newcommand{\PGC}{\mathrm{PGC}}
\newcommand{\GC}{\mathrm{GC}}
\newcommand{\U}{\mathrm{U}}
\newcommand{\CC}{\mathrm{C}}
\newcommand{\Bb}{\mathrm{B}}
\newcommand{\K}{\mathrm{K}}

\newcommand{\mtx}[2]
{\left(\!\!\!
	\begin{array}{cc}
		#1  \\
		#2 
	\end{array}
	\!\!\!\right)}

\newcommand {\mmtx}[9]
{\left(
	\begin{array}{ccc}
		#1 & #2 & #3 \\
		#4 & #5 & #6 \\
		#7 & #8 & #9
	\end{array}
	\right)}

\newcommand {\mtxx}[4]
{\left(\!
	\begin{array}{cc}
		\!\!#1 & \!\!#2 \\
		\!\!#3 & \!\!#4
	\end{array}\!\!
	\right)}

\newtheoremstyle{athm}
{}
{}
{\itshape}
{}
{\scshape}
{}
{.5em}
{\thmnote{#3}}
\theoremstyle{athm}

\begin{document}
	\title[Homology of projective linear group $\PGL_2$]
	{The low dimensional homology of projective linear group of degree two}
	
	\author{Behrooz Mirzaii}
	\author{Elvis Torres P\'erez}
	
	\address{\sf 
		Instituto de Ci\^encias Matem\'aticas e de Computa\c{c}\~ao (ICMC), Universidade de S\~ao Paulo, 
		S\~ao Carlos, Brasil}
	\email{bmirzaii@icmc.usp.br}
	\email{elvistorres@usp.br}
	
	\begin{abstract}
		In this article we study the low dimensional homology of the projective linear group $\PGL_2(A)$ over a 
		commutative ring $A$. In particular, we prove a Bloch-Wigner type exact sequence over local domains. As 
		application we prove that $H_2(\PGL_2(A),\z\half)\simeq \K_2(A)\half$ and $H_3(\PGL_2(A),\z\half)\simeq 
		\K_3^\ind(A)\half$, provided $|A/\mmm_A|\neq 2,3,4, 8$.
	\end{abstract}
	\maketitle
	
	Let $A$ be a commutative ring with 1. Let $\GE_2(A)$ be the subgroup of $\GL_2(A)$ generated by
	elementary and diagonal matrices. We say that $A$ is a $\GE_2$-ring if $\GE_2(A)=\GL_2(A)$. 
	This is equivalent to the condition that $\Ee_2(A)=\SL_2(A)$.
	
	A ring $A$ is called universal for $\GE_2$ if the unstable $K$-group $\K_2(2,A)$ is generated by 
	Steinberg symbols (see Section \ref{sec2}).  We say that $A$ is an universal $\GE_2$-ring if it 
	is a $\GE_2$-ring and is universal for $\GE_2$.
	If $G$ is any subgroup of $\GL_2(A$) containing the central subgroup $Z=\aa I_2$ of scalar
	matrices, then we will let ${\rm PG}$ denote the quotient group $G/Z$.
	
	As our first main result we show that for any commutative ring $A$, we have the exact sequence
	\begin{equation}\label{first}
		H_2(\PGE_2(A),\z) \arr \bigg(\displaystyle{\frac{\K_2(2,A)}{\CC(2, A)}}\bigg)_{\PGE_2(A)}^\ab
		\arr A_\aa \arr H_1(\PGE_2(A),\z) \arr \GG_A \arr 1,
	\end{equation}
	where $\CC(2,A)$ is the central subgroup of $\K_2(2,A)$ generated by Steinberg symbols, $\GG_A$ is 
	the square class group of $A$, i.e. $\GG_A:=\aa/(\aa)^2$, and $A_\aa:=A/\lan a-1:a\in \aa\ran$
	(see Theorem \ref{E1S}). It follows from this that if $A$ is an universal $\GE_2$-ring, then 
	\[
	H_1(\PGL_2(A),\z)\simeq \GG_A\oplus A_\aa.
	\]
	
	As our second main result we show that if $A$ is an universal $\GE_2$-ring, then we have
	the exact sequence
	\begin{equation}\label{second}
		H_3(\PGL_2(A),\z) \arr \PP(A) \overset{\lambda}{\larr} H_2(\PB_2(A),\z) \arr H_2(\PGL_2(A),\z) 
		\arr \mu_2(A) \arr 1,
	\end{equation}
	where $\PP(A)$ is the scissors congruence group of $A$ and $\PB_2(A)$ is the group of upper 
	triangular matrices in $\PGL_2(A)$ (for the general statement see Theorem~\ref{ES} and 
	Corollary~\ref{H1-uni}). 
	
	Let $A$ be a local ring such that $|A/\mmm_A|\neq 2,3,4$. Then $H_2(\PB_2(A),\z)\simeq \aa\wedge\aa$
	(Proposition~\ref{B2-T2}) and we show that the map $\lambda$ is given by 
	\[
	\lambda([a])=2(a\wedge (1-a))
	\]
	(Proposition \ref{d21}). As an application we show that if $A$ is a local domain (local ring) such 
	that $|A/\mmm_A|\neq 2, 3, 4$ ($|A/\mmm_A|\neq 2, 3, 4, 5,8,9,16$), then
	\begin{equation}\label{K2}
		\begin{array}{c}
			H_2(\PGL_2(A),\z\half)\simeq \K_2(A)\half.
		\end{array}
	\end{equation}
	
	Let $\BB_E(A)$ be the kernel of $\lambda$. Then as our third main result we show that if 
	$A$ is a local domain such that $|A/\mmm_A|\neq 2,3,4,8$, then we obtain the sequence
	\begin{equation}\label{PBW}
		\begin{array}{c}
			0\arr \tors(\mu(A),\mu(A)) \arr H_3(\PGL_2(A),\z) \arr \BB_E(A) \arr 0,
		\end{array}
	\end{equation}
	which is exact at every term except possibly at the term $H_3(\PGL_2(A),\z)$, where the homology of the 
	sequence is annihilated by $4$ (see Theorem \ref{main-thm} for the general statement).  As an application 
	we prove the Bloch-Wigner exact sequence
	\begin{equation}\label{P-BW}
		\begin{array}{c}
			0\arr \tors(\mu(A),\mu(A))\half \arr H_3(\PGL_2(A),\z\half) \arr \BB(A)\half \arr 0,
		\end{array}
	\end{equation}
	where $\BB(A)\se \PP(A)$ is the Bloch group of $A$. As an application of this exact sequence we show that
	\[
	\begin{array}{c}
		H_3(\PGL_2(A),\z\half)\simeq \K_3^\ind(A)\half.
	\end{array}
	\]
	
	The earliest version of the celebrated Bloch-Wigner exact sequence that we found in the literature is 
	the exact sequence
	\[
	\begin{array}{c}
		0 \arr \tors(\mu(\C),\mu(\C)) \arr H_3(\PGL_2(\C),\z) \arr \BB(\C) \arr 0
	\end{array}
	\]
	(see \cite[Theorem 4.10]{dupont-sah1982}). The exact sequence (\ref{P-BW}) can be seen as a generalization 
	of this classical result to local domains. As we will see, in general the coefficient $\z\half$ can not be 
	replaced with integral coefficients $\z$, even over infinite fields (see for example Proposition \ref{localfield}). 
	Moreover, we study the sequence (\ref{PBW}) over quadratically closed fields, real closed fields, finite fields 
	and non-dyadic local fields (see Propositions \ref{classical}, \ref{localfield}). Finally we prove a Bloch-Wigner 
	type exact sequence for $\PGL_2(\z)$ and $\PGL_2(\z\half)$.
	
	Here we outline the organization of the present paper. In Section~\ref{sec1}, we  recall some needed results from 
	the literature over algebraic $K$-groups, the scissors congruence group and the Bloch-Wigner exact sequence. 
	In Section~\ref{sec2} we recall the Steinberg group $\St(2,A)$, the $K$-group $\K_2(2,A)$ and give some of its
	basic properties. In Section~\ref{sec3} we give a detailed account of the action of $\PGE_2(A)$ over $\K_2(2,A)$,
	construct the important map 
	\[
	\kappa: \bigg(\displaystyle{\frac{\K_2(2,A)}{\CC(2, A)}}\bigg)_{\PGE_2(A)}^\ab \arr A_\aa
	\]
	and prove our first main result, i.e. the exactness of the sequence (\ref{first}). In Section~\ref{sec4}, we study two 
	chain complexes $Y_\bullet(A^2)\subseteq L_\bullet(A^2)$ made out of unimodular vectors in $A^2$ which are columns
	of matrices in $\GE_2(A)$ and $\GL_2(A)$, respectively, and study the connection between their homology groups. 
	In Section \ref{sec5} we study the connection between the first homology group of these complexes and the group 
	$\bigg(\displaystyle{\frac{\K_2(2,A)}{\CC(2, A)}}\bigg)_{\PGE_2(A)}^\ab$. In Section~\ref{sec6}, we introduce 
	and study a spectral sequence which will be our main tool in handling the second and the third homology groups 
	of $\PGE_2(A)$. In Section~\ref{sec7}, we study certain terms of the spectral sequence 
	and prove the exactness of the sequence (\ref{second}). In Section~\ref{sec8}, the homology groups of $\PB_2(A)$ 
	have been studied. In Section~\ref{sec9}, we calculate the map $\lambda$ and prove the isomorphism (\ref{K2}). In 
	Section~\ref{sec10} we prove our claim about the sequence (\ref{PBW}) and  present the proof of the Bloch-Wigner exact 
	sequence (\ref{P-BW}). Moreover, we prove a Bloch-Wigner type exact sequence over finite fields, real closed 
	fields, non-dyadic local fields and the euclidean domains $\z$ and $\z\half$.
	\\
	~\\
	{\bf Notations.}
	In this paper all rings are commutative, except possibly group rings, and have the unit element $1$. For a 
	commutative ring $A$ let $\GL_2(A)$ be the group of invertible matrices of degree two. If ${\rm G}(A)$ is a 
	subgroup of $\GL_2(A)$ which contains $\aa I_2=Z(\GL_2(A))$, by ${\rm PG}(A)$ we mean ${\rm G}(A)/\aa I_2$. 
	Let $\mu(A)$ denote the group of roots of unity in $A$, i.e.
	\[
	\mu(A):=\{a\in A: \text{there is $n\in \N$ such that $a^n=1$} \},
	\]
	and $\mu_2(A):=\{a\in A: a^2=1\}$. Let $\GG_A:=\aa/(\aa)^2$. The element of $\GG_A$ represented by 
	$a \in \aa$ is denoted by $\lan a \ran$. If $\BB \arr \Aa$ is a homomorphism of abelian groups, by $\Aa/\BB$ 
	we mean $\coker(\BB \arr \Aa)$. For an abelian group $\Aa$, by $\Aa\half$ we mean $\Aa\otimes_\z \z\half$.\\
	~\\
	{\bf Acknowledgements.}
	We would like to thank the anonymous referee for many valuable comments and suggestions which improved 
	some of our main results. Moreover his suggestions greatly improved the presentation of the present 
	article. The detailed action of $\PGE_2(A)$ on $\K_2(2,A)$, the description of the map $\kappa$ 
	in Section \ref{sec3} and the use of the functor $\Ind_{\PGE_2}^{\PGL_2}$ in connection with the 
	complexes $L_\bullet(A^2)$ and $Y_\bullet(A^2)$ in Sections \ref{sec4} and \ref{sec5} were suggested by the referee. 
	
	\section{Algebraic {\it K}-theory and scissors congruence group}\label{sec1}
	
	Let $A$ be a commutative ring. For any non-negative integer $n\geq 1$, we associate two 
	type $K$-groups to $A$: Quillen's $K$-group $\K_n(A)$ and Milnor's $K$-group $\K_n^M(A)$. 
	
	Quillen's $K$-group $\K_n(A)$ is defined as the $n$-th homotopy group of the plus-construction 
	of the classifying space of the stable linear group $\GL(A)$, with respect to the perfect 
	elementary subgroup $\Ee(A)$:
	\[
	\K_n(A):=\pi_n(B\GL(A)^+).
	\]
	Since $B\Ee(A)^+$ is homotopy equivalent to the universal cover of $B\GL(A)^+$, for $n\geq 2$ we have
	\[
	\K_n(A)\simeq\pi_n(B\Ee(A)^+).
	\]
	The Hurewicz map in algebraic topology induces the commutative diagram (for $n\geq 2)$
	\[
	\begin{tikzcd}
		& & H_n(\Ee(A),\z) \ar[dd, hook]\\
		&\K_n(A)  \arrow[rd, "h_n"] \ar[ru,"h_n'"] && \\
		& & H_n(\GL(A),\z)
	\end{tikzcd}
	\]
	If $A'$ is another commutative ring, there is a natural anti-commutative product map
	\[
	\K_m(A) \otimes_\z \K_n(A') \arr \K_{m+n}(A\otimes_\z A'),  \ \ \ \ \ x\otimes y \mapsto x\star y.
	\]
	When $A'=A$ and $\eta:A\otimes_\z A \arr A$ is given by $a\otimes b \mapsto ab$, then we have 
	the product map
	\[
	\K_m(A) \otimes_\z \K_n(A) \overset{\eta_\ast\circ\star}{\larr} 
	\K_{n+m}(A), \ \ \ \ \ x\otimes y \mapsto \eta_\ast(x\star y).
	\]
	For more on these $K$-groups and the construction of the product map see \cite[Chap. 2]{srinivas1996}.
	
	The $n$-th Milnor $K$-group $\K_n^M(A)$ is defined as the abelian group generated by symbols 
	$\{a_1,\dots, a_n\}$, $a_i \in \aa$, subject to the following relations\\
	
	\par (i) $\{a_1,\dots, a_ia_i', \dots a_n\} = \{a_1, \dots, a_i,\dots, a_n\}+ 
	\{a_1, \dots, a_i',\dots, a_n\}$, for any $1\leq i\leq n$,
	\par (ii) $\{a_1,\dots,a_n\}=0$ if there exist $i, j$, $i\neq j$, such that $a_i+a_j=0$ or $1$.\\
	~\\
	Clearly we have the anti-commutative product map
	\[
	\K_m^M(A) \otimes_\z \K_n^M(A) \arr \K_{m+n}^M(A),
	\]
	\[
	\{a_1,\dots, a_m\}\otimes \{b_1,\dots, b_n\}\mapsto\{a_1,\dots, a_m, b_1,\dots, b_n\}.
	\]
	It can be shown that 
	\[
	\K_1(A) \overset{h_1}{\simeq} H_1(\GL(A),\z)\simeq \GL(A)/\Ee(A), \ \ \ \ 
	\K_2(A)\overset{h_2'}{\simeq} H_2(\Ee(A),\z).
	\]
	For $n=1$, we have the natural homomorphism
	\[
	\K_1^M(A) \arr \K_1(A), \ \ \ \ \ \{a\} \mapsto {\mtxx{a}{0}{0}{1}}.
	\]
	The determinant induces the isomorphism 
	\[
	\K_1(A)\simeq \aa \times {\rm SK}_1(A)\simeq \K_1^M(A) \times {\rm SK}_1(A),
	\]
	where ${\rm SK}_1(A):= \SL(A)/\Ee(A)$. 
	If $A$ is a local ring, then ${\rm SK}_1(A)=1$ and thus 
	\[
	\K_1(A)\simeq K_1^M(A).
	\]
	For $n=2$ we have the natural homomorphism 
	\[
	\K_2^M(A) \arr \K_2(A),  \ \ \ \ \{ a,b\} \mapsto 
	\eta_\ast\bigg({\mmtx{a}{0}{0}{0}{a^{-1}}{0}{0}{0}{1}}\star{\mmtx{b}{0}{0} {0}{1}{0} {0}{0}{b^{-1}}}\bigg).
	\]
	The following result is well-known.
	
	\begin{thm}[Matsumoto-van der Kallen]\label{MV}
		Let $A$ be either a field or a local ring such that its residue field has more than five elements. The natural 
		homomorphism $\K_2^M(A)\arr \K_2(A)$ is an isomorphism
		\[
		(\aa \otimes_\z \aa)/\lan a\otimes (1-a): a(1-a)\in \aa \ran \simeq \K_2^M(A)\simeq \K_2(A).
		\]
	\end{thm}
	\begin{proof}
		See \cite[Theorem 1.14]{srinivas1996} and \cite[Proposition 3.2]{mirzaii2017}.
	\end{proof}
	
	Using products of $K$-groups, one can show that for any positive integer $n$, there is a natural map
	\[
	\psi_n: \K_n^M(A) \arr \K_n(A).
	\]
	For the group $\K_3(A)$ we have the following general result.
	
	\begin{thm}[Suslin]\label{sus}
		For any ring $A$ we have the exact sequence
		\[
		\K_1(\z)\otimes_\z \K_2(A) \overset{\star}{\larr} \K_3(A) \overset{h_3'}{\larr} H_3(\Ee(A),\z) \arr 0.
		\]
	\end{thm}
	\begin{proof}
		See \cite[Corollary 5.2]{suslin1991} 
	\end{proof}

	Let $\WW_A:=\{a\in A: a(1-a)\in \aa\}$. By definition, the {\it classical scissors congruence group} 
	$\PP(A)$ of $A$ is the quotient of the free abelian group generated by symbols $[a]$, $a\in \WW_A$, 
	by the subgroup generated by the elements
	\[
	[a] -[b]+\bigg[\frac{b}{a}\bigg]-\bigg[\frac{1- a^{-1}}{1- b^{-1}}\bigg]+ \bigg[\frac{1-a}{1-b}\bigg],
	\]
	where $a, b, a/b  \in \WW_A$. Let 
	\[
	S_\z^2(\aa):=(\aa\otimes_\z \aa)/\lan a\otimes b+b\otimes a:a,b \in \aa\ran.
	\]
	The map
	\[
	\lambda: \PP(A) \arr S_\z^2(\aa), \ \ \ \ \ [a]\mapsto a\otimes (1-a)
	\]
	is well-defined. The kernel of $\lambda$ is called the {\it Bloch group} of $A$ and is denoted by 
	$\BB(A)$. If $A$ is either a field or a local ring such that its residue field has more than five 
	elements, then we have the exact sequence
	\[
	0 \arr \BB(A) \arr \PP(A) \arr S_\z^2(\aa) \arr \K_2^M(A) \arr 0.
	\]
	
	The group $\K_3(A)$ is closely related to the Bloch group of $A$. Over a local ring, the indecomposable 
	part of $\K_3(A)$ is defined as follows:
	\[
	\K_3^\ind(A):=\K_3(A)/\K_3^M(A).
	\]
	
	\begin{thm}[A Bloch-Wigner exact sequence]
		Let $A$ be either a field or a local domain such that its residue field has more than $9$ elements. 
		Then there is a natural exact sequence
		\[
		0 \arr \tors(\mu(A),\mu(A))^\sim \arr \K_3^\ind(A) \arr \BB(A) \arr 0, 
		\]
		where $\tors(\mu(A),\mu(A))^\sim$ is the unique nontrivial extension of $\tors(\mu(A),\mu(A))$ by $\mu_2(A)$.
	\end{thm}
	\begin{proof}
		The case of infinite fields has been proved by Suslin in \cite[Theorem 5.2]{suslin1991} and the 
		case of finite fields has been settled by Hutchinson in \cite[Corollary 7.5]{hutchinson-2013}. 
		The case of local rings has been dealt in \cite[Theorem 6.1]{mirzaii2017}.    
	\end{proof}
	
	Let $A$ be either a field or a local domain such that its residue field has more than five 
	elements. Then by Theorem \ref{MV}, $\K_2(A)\simeq \K_2^M(A)$. Since $\K_1(\z)\simeq\{\pm 1\}$ 
	\cite[Example 1.9(vii)]{srinivas1996}, we have
	\[
	\im(\K_1(\z)\otimes_\z \K_2(A) \overset{\star}{\larr} \K_3(A)) \se \im(\K_3^M(A) \arr \K_3(A)).
	\]
	Let $\alpha_A$ be the following composite
	\[
	H_3(\SL_2(A),\z)_\aa \arr H_3(\SL(A),\z) \simeq \K_3(A)/(\K_1(\z)\otimes_\z \K_2(A)) \arr \K_3^\ind(A).
	\]
	Note that over local rings $\Ee(A)=\SL(A)$ and $\aa\simeq \GL(A)/\Ee(A)$ acts trivially on the 
	group $H_3(\SL(A),\z)$. The following question was asked by Suslin (see \cite[Question 4.4]{sah1989}).
	
	\begin{question}
		For an infinite field $F$, is the map $\alpha_F :H_3(\SL_2(F),\z)_\ff \arr \K_3^\ind(F)$ an isomorphism?
	\end{question}
	
	Hutchinson and Tao proved that $\alpha_F$ always is surjective \cite[Lemma 5.1]{hutchinson-tao2009}. 
	The answer of the above question is true for all finite fields except for $F=\F_2,\F_3,\F_4, \F_8$ 
	(see \cite[Proposition 6.4, Example 6.6]{mirzaii2017}). For more on the above question see \cite{mirzaii2015}.
	
	\begin{thm}\label{kind-sl2}
		Let $A$ be a local domain such that $|A/\mmm_A|\neq 2,3,4,8$. Then the map 
		\[
		\begin{array}{c}
			\alpha_A:H_3(\SL_2(A),\z\half)_\aa \arr  \K_3^\ind(A)\half
		\end{array}
		\]
		is an isomorphism.
	\end{thm}
	\begin{proof}
		See \cite[Theorem 3.7]{mirzaii2012}, \cite[Theorem 5.4]{mirzaii2017} and \cite[Theorem 6.4]{mirzaii-2008}.
	\end{proof}
	
	\section{ Elementary matrices and the Steinberg group of degree two
	}\label{sec2}
	
	Let $A$ be a commutative ring. The elementary group of degree two over $A$, denoted by $\Ee_2(A)$, is the subgroup of $\GL_2(A)$ 
	generated by the elementary matrices
	\[
	E_{12}(a):=\begin{pmatrix}
		1 & a\\
		0 & 1
	\end{pmatrix}, \ \ \ \ \ 
	E_{21}(a):=\begin{pmatrix}
		1 & 0\\
		a & 1
	\end{pmatrix}, \ \ \ \ \ a\in A.
	\]
	The elementary matrices satisfy the following relations
	\begin{itemize}
		\item[(a)] $E_{ij}(x)E_{ij}(y)=E_{ij}(x+y)$ for any $x,y \in A$,
		\item[(b)] $W_{ij}(u)E_{ji}(x)W_{ij}(u)^{-1}=E_{ij}(-u^2x)$, for any $u\in \aa$ and $x\in A$,\\
		where $W_{ij}(u):=E_{ij}(u)E_{ji}(-u^{-1})E_{ij}(u)$.
	\end{itemize}
	The {\it Steinberg group} of $A$, denoted by $\St(2,A)$, is the group with generators $x_{12}(r)$ and 
	$x_{21}(s)$, $r,s\in A$, subject to the relations
	\begin{itemize}
		\item[($\alpha$)] $x_{ij}(r)x_{ij}(s)=x_{ij}(r+s)$ for any $r,s \in A$,
		\item[($\beta$)] $w_{ij}(u)x_{ji}(r)w_{ij}(u)^{-1}=x_{ij}(-u^2r)$, for any $u\in \aa$ and $r\in A$,\\
		where $w_{ij}(u):=x_{ij}(u)x_{ji}(-u^{-1})x_{ij}(u)$.
	\end{itemize}
	The natural map
	\[
	\phi: \St(2,A) \arr \Ee_2(A), \ \ \ \ \ x_{ij}(r) \mapsto E_{ij}(r)
	\]
	is a well defined epimorphism. The kernel of this map is denoted by $\K_2(2,A)$:
	\[
	\K_2(2,A):=\ker(\phi).
	\]
	Always there is a natural map
	\[
	\K_2(2,A) \arr \K_2(A),
	\]
	which in general neither is surjective nor injective. If $A$ is a local ring, then this map always is 
	surjectve \cite[Theorem 2.13]{stein1973}.
	
	For any $u\in \aa$, let
	\[
	h_{ij}(u):=w_{ij}(u)w_{ij}(-1).
	\]
	It is not difficult to see that $h_{ij}(u)^{-1}=h_{ji}(u)$ \cite[Corollary A.5]{hutchinson2022}. For any 
	$u,v\in \aa$, the element
	\[
	\{u,v\}_{ij}:=h_{ij}(uv)h_{ij}(u)^{-1}h_{ij}(v)^{-1} 
	\]
	lies in $\K_2(2, A)$ and in the center of $\St(2, A)$  \cite[\S 9]{D-S1973}. It is straightforward to check that
	$\{u,v\}_{ji}=\{v,u\}_{ij}^{-1}$.
	An element of form 
	\[
	\{v,u\}:=\{v,u\}_{12}=h_{12}(uv)h_{12}^{-1}(u)h_{12}(v)^{-1}
	\]
	is called a {\it Steinberg symbol} in $\K_2(2, A)$.
	
	Let $\CC(2,A)$ be the subgroup of $\K_2(2, A)$ generated by the Steinberg 
	symbols $\{u,v\}$, $u,v\in \aa$. Then $\CC(2, A)$ is a central subgroup of $\K_2(2, A)$.
	
	We say that $A$ is {\it universal for} $\GE_2$ if $\K_2(2, A)=\CC(2, A)$. This definition of universal for $\GE_2$ 
	is equivalent to the original definition of Cohn in \cite[page 8]{cohn1966}. For a proof of this fact see 
	\cite[Appendix A]{hutchinson2022}. A commutative semilocal ring is universal for $\GE_2$ if and only if none of the rings  
	$\z/2 \times \z/2$ and $\z/6$ is a direct factor of $A/J(A)$, where $J(A)$ is the  Jacobson radical of $A$ 
	\cite[Theorem 2.14]{menal1979}.
	
	The elementary group $\Ee_2(A)$ is generated by the matrices 
	\[
	E(a):={\mtxx a 1 {-1} 0}, \ \ \ \ \ \ a \in  A.
	\]
	In fact
	\[
	E_{12}(a)=E(-a)E(0)^{-1},\ \ \ \ \ \ E_{21}(a)=E(0)^{-1} E(a),
	\]
	\[
	E(0)=E_{12}(1)E_{21}(-1)E_{12}(1).
	\]
	For any $a\in \aa$, let 
	\[
	D(a):=\begin{pmatrix}
		a & 0\\
		0 & a^{-1}
	\end{pmatrix}\in {\rm D}_2(A).
	\]
	Since $D(-a)=E(a)E(a^{-1})E(a)$, we have $D(a)\in \Ee_2(A)$. It is straightforward to check that
	\begin{itemize}
		\item[(1)] $E(x)E(0)E(y)=D(-1)E(x+y)$,
		\item[(2)] $E(x)D(a)=D(a^{-1})E(a^2x)$,
		\item[(3)] $D(ab)D(a^{-1})D(b^{-1})=1$,
	\end{itemize}
	where $x,y\in A$ and $a, b\in \aa$. Let $\CC(A)$ be the group generated by symbols $\ep(a)$, 
	$a\in A$, subject to the relations
	\begin{itemize}
		\item[(i)] $\ep(x)\ep(0)\ep(y)=h(-1)\ep(x+y)$ for any $x,y \in A$,
		\item[(ii)] $\ep(x)h(a)=h(a^{-1})\ep(a^2x)$, for any $x\in A$ and $a\in \aa$,
		\item[(iii)] $h(ab)h(a^{-1})h(b^{-1})\!=\!1$ for any $a, b\in\aa$,
	\end{itemize}
	where 
	\[
	h(a):=\ep(-a)\ep(-a^{-1})\ep(-a).
	\]
	Note that by (iii), $h(1)=1$ and $h(-1)^2=1$. Moreover 
	$\ep(-1)^3=h(1)=1$ and $\ep(1)^3=h(-1)$. There is a natural surjective map 
	\[
	{\rm C}(A) \arr \Ee_2(A),  \ \ \ \ \ep(x)\mapsto E(x).
	\]
	We denote the kernel of this map by $\U(A)$. Thus we have the extension
	\[
	1 \arr {\rm U}(A) \arr {\rm C}(A) \arr \Ee_2(A) \arr 1.
	\]
	
	\begin{prp}[Hutchinson]\label{hut}
		Let $A$ be a commutative ring. Then the map $\St(2,A) \arr \CC(A)$ given by $x_{12}(a)\mapsto \varepsilon(-a)\varepsilon(0)^3$ 
		and $x_{21}(a)\mapsto\varepsilon(0)^3\varepsilon(a)$ induces isomorphisms 
		\[
		\displaystyle \frac{\St(2,A)}{\CC(2, A)}\simeq \CC(A),\ \ \ \ \displaystyle \frac{\K_2(2,A)}{\CC(2, A)}\simeq \U(A).
		\]
	\end{prp}
	\begin{proof}
		See \cite[Theorem A.14, App. A]{hutchinson2022}.    
	\end{proof}
	
	It follows from this theorem that $A$ is universal for $\GE_2$ if and only if $\U(A)=1$.

	\section{ The group \texorpdfstring{$\K_2(2,A)$}{Lg} and the abelianization of \texorpdfstring{$\GE_2(A)$}{Lg}}\label{sec3}
	
	Let $A$ be a commutative ring. 
	Let ${\rm D}_2(A)$ be the subgroup of $\GL_2(A)$ generated by diagonal matrices and 
	let $\GE_2(A)$ be the subgroup of $\GL_2(A)$ generated by ${\rm D}_2(A)$ and $\Ee_2(A)$. 
	It is easy to see that $\Ee_2(A)$ is normal in $\GE_2(A)$ \cite[Proposition 2.1]{cohn1966} 
	and the center of $\GE_2(A)$ is $\aa I_2$. Observe that we have the split extensions
	\[
	1 \arr \SL_2(A) \arr \GL_2(A) \overset{\det}{\arr} \aa \arr  1 , \ \ \ \ \ 
	1 \arr \Ee_2(A) \arr \GE_2(A) \overset{\det}{\arr} \aa \arr  1,
	\]
	and thus
	\[
	\GL_2(A)=\SL_2(A)  \rtimes d(\aa), \ \ \ \ \GE_2(A)= \Ee_2(A)  \rtimes d(\aa),
	\]
	where 
	\[
	d(\aa)=\{d(a):=\diag(a,1):a\in \aa\}\simeq \aa.
	\]
	We say that $A$ is a $\GE_2$-{\it ring} if $\GE_2(A)=\GL_2(A)$ (or equivalently $\Ee_2(A)=\SL_2(A)$). 
	Semilocal rings and Euclidean domains are $\GE_2$-rings \cite[p. 245]{silvester1982}, \cite[\S 2]{cohn1966}. 
	
	A ring $A$ is called an {\it universal $\GE_2$-ring} if it is a $\GE_2$-ring and is universal for $\GE_2$.
	A semilocal ring is an universal $\GE_2$-ring if none of the rings $\z/2 \times \z/2$ and $\z/6$ is a direct 
	factor of $A/J(A)$. In particular, any local ring is an universal $\GE_2$-ring.
	For more example of $\GE_2$-rings and rings universal for $\GE_2$ see \cite{cohn1966}, \cite{hutchinson2022} and 
	\cite{B-E2024}.
	
	There is a natural action of $\PGE_2(A)$ on $\K_2(2, A)$. Here we give a detailed description
	of this action. From the extension
	\[
	1 \arr \K_2(2, A) \arr \St(2, A) \arr \Ee_2(A) \arr 1,
	\]
	we see that $\Ee_2(A)$ acts naturally on  $\K_2(2, A)$. More explicitly, $E_{12}(t)$ acts as 
	conjugation by $x_{12}(t)$ and $E_{21}(t)$ acts as conjugation by $x_{21}(t)$. Note that $D(a)=\diag(a, a^{-1})$ 
	acts as conjugation by $h_{12}(a)$. It is straightforward to check that
	\[
	x_{12}(t)^{h_{12}(a)}=x_{12}(a^{-2}t), \ \ \ \ x_{21}(t)^{h_{12}(a)}=x_{21}(a^2t).
	\]
	In particular, the scalar matrix $-I_2=D(-1) \in \Ee_2(A)$ acts trivially 
	on $\K_2(2, A)$.
	
	For $a \in \aa$, let $d(a) := \diag(a, 1) \in \GE_2(A)$. For any $t \in A$, 
	\[
	E_{12}(t)^{d(a)}=E_{12}(a^{-1}t), \ \ \ \ E_{21}(t)^{d(a)} = E_{21}(at). 
	\]
	It is straightforward to verify that there is a compatible well-defined action of 
	$d(\aa)$ on $\St(2, A)$ determined by
	\[
	x_{12}(t)^{d(a)}:= x_{12}(a^{-1}t),\ \ \ \ \ x_{21}(t)^{d(a)}:= x_{21}(at).
	\]
	One verifies easily that this proposed action preserves the two defining families of relations of $\St(2, A)$.
	
	This implies that $\GE_2(A)=\Ee_2(A) \rtimes d(\aa)$ acts on $\St(2,A)$ via the conjugation formula given above. 
	Using this action we define 
	\[
	\GSt(2, A) := \St(2, A) \rtimes d(\aa)
	\]
	and extend the canonical epimorphism $\phi : \St(2, A) \arr \Ee_2(A)$ to
	a surjective homomorphism 
	\[
	\Phi : \GSt(2, A) \arr \GE_2(A).
	\]
	Furthermore, the inclusion $\St(2, A) \arr \GSt(2, A)$ induces an isomorphism 
	\[
	\K_2(2, A)=\ker\phi \simeq \ker\Phi.
	\]
	Thus we have the extension 
	\[
	1 \arr \K_2(2, A) \arr \GSt(2, A) \arr \GE_2(A) \arr 1.
	\]
	Thus $\GE_2(A)$ acts by conjugation on $\K_2(2, A)$ which is compatible with the above action of $\Ee_2(A)$.
	
	The matrix $d(a)$ acts by the formula given above. If we let $d'(a) :=\diag(1, a)$, then 
	$d'(a) = \diag(a^{-1}, a)d(a)$ which lifts to $(h_{12}(a^{-1}),d(a))$ in $\GSt(2, A)$. Then 
	this matrix acts on $\St(2,A)$ via the formulas
	\[
	x_{12}(t)^{d'(a)}= x_{12}(at), \ \ \ \  x_{21}(t)^{d'(a)}= x_{21}(a^{-1}t).
	\]
	It follows in turn that the scalar matrices $aI_2 = d(a)d'(a)$ act trivially.
	Hence the above action descends to an action of $\PGE_2(A)$ on $\K_2(2, A)$.
	
	Since $\CC(2, A)$ is central in $\St(2, A)$, the action by conjugation of $\St(2, A)$ on $\CC(2, A)$ is trivial and 
	hence $\Ee_2(A)$ acts trivially on the image, $\overline{\CC}(2, A)$, of $\CC(2, A)$ in $\K_2(2, A)^\ab$. It can be 
	easily verified that the action of $d(\aa)$ on $\K_2(2, A)$ induces an action on $\CC(2, A)$ given by 
	\[
	\{u, v\}^{d(a)} = \{u, a^{-1}\}^{-1}\{u, a^{-1}v\}.
	\]
	Hence the action of $\PGE_2(A)$ on $\K_2(2, A)$ induces an action on the group
	\[
	\displaystyle \frac{\K_2(2,A)^\ab}{\overline{\CC}(2, A)}\simeq\bigg(\frac{\K_2(2,A)}{\CC(2,A)}\bigg)^\ab.
	\]
	
	Let $A_\aa:=A/\lan a-1:a\in \aa\ran$. Since $-2=(-1)-1$, we have $2=0$ in $A_\aa$. Hence we have a natural map
	\[
	A/2A \arr A_\aa.
	\]
	
	There is a natural homomorphism 
	\[
	f : \St(2, A) \arr A_\aa
	\]
	which sends both $x_{12}(t)$ and $x_{21}(t)$ to (the class of) $t$, for $t \in A$. Observe that this
	map sends the elements $w_{ij} (a)$ to $1$. Thus it sends $h_{ij}(a)$ to $2 = 0$.
	Therefore $f$ extends to a homomorphism 
	\[
	f : \GSt(2, A) = \St(2, A)\rtimes d(\aa) \arr A_\aa, \ \ \ \ (x, d) \mapsto f(x).
	\]
	For example $x_{12}(t)^{d(a)} = x_{12}(a^{-1}t)$ and both sides map to $t$. Now the diagram
	\[
	\begin{tikzcd}
		& \GSt(2, A) \ar[dl, "f"] \ar[dr] & \\
		A_\aa \ar[rr] & & \PGE_2(A)^\ab
	\end{tikzcd}
	\]
	commutes: Here the map $A_\aa \arr \PGE_2(A)^\ab$ sends $t$ to $E_{12}(t)$ and observe  that in $\PGE_2(A)^\ab$
	we have
	\[
	E_{12}(t) = E_{12}(t)^{W_{12}(1)} = E_{21}(-t) = E_{21}(-t)^{d(-1)} = E_{21}(t).
	\]
	We denote the restriction of $f$ to $\K_2(2, A)$ again by $f$:
	\[
	f:\K_2(2, A) \arr  A_\aa.
	\]
	Note that since $f(h_{ij}(u)) = 0$ for all $u\in \aa$, then $f(\{u, v\}) = 0$ for all 
	$u, v \in \aa$. Thus $f$ naturally defines a map
	\[
	\kappa=\overline{f}: \bigg(\frac{\K_2(2, A)}{\CC(2, A)}\bigg)^\ab \arr A_\aa.
	\]
	
	Using the action of $d(\aa)$ on $\St(2,A)$, we can define an action of $d(\aa)$ on $\CC(A)$. More precisely,
	for any $a \in \aa$ and $x \in A$ we have
	\[
	E(x)^{d(a)}=D(a)E(a^{-1}x).
	\]
	Thus the compatible action of $d(\aa) \simeq \aa$ on $\CC(A)$ is determined by
	\[
	\ep(x)^{d(a)}=h(a)\ep(a^{-1}x).
	\]
	It is easy to verify that $h(b)^{d(a)}=h(b)$. If $\GC(A):=\CC(A)\rtimes d(\aa)$, then we have the extension
	\[
	1 \arr \U(A) \arr \GC(A) \arr \GE_2(A) \arr 1.
	\]
	Observe that we have the natural morphism of extensions
	\[
	\begin{tikzcd}
		1 \ar[r] & \K_2(2,A) \ar[r]\ar[d]& \GSt(2,A) \ar[r]\ar[d] & \GE_2(A) \ar[r]\ar[d]& 1\\
		1 \ar[r] & \U(A) \ar[r] & \GC(A) \ar[r] & \GE_2(A) \ar[r]&  1.
	\end{tikzcd}
	\]
	It is easy to check that the action of $aI_2$ on $\GC(A)$ is trivial. Note that $aI_2=D(a^{-1})d(a^2)$.
	
	Let ${\rm L}(A)$ be the subset of $\GC(A)$ consisting of the elements $(h(a^{-1}), d(a^2))$,  $a\in \aa$.
	Note that $(h(a^{-1}), d(a^2))$ lies in the center of $\GC(A)$:
	\begin{align*}
		(h(a^{-1}), d(a^2))(\varepsilon(x), d(b))&=(h(a^{-1})\varepsilon(x)^{d(a^2)}, d(a^2b))\\
		&=(h(a^{-1})h(a^2)\varepsilon(a^{-2}x), d(a^2b))\\
		&=(h(a)h(a^{-1})\varepsilon(x)h(a^{-1}), d(a^2b))\\
		&=(\varepsilon(x)h(a^{-1}), d(a^2b))\\
		&=(\varepsilon(x)h(a^{-1})^{d(b)}, d(ba^2))\\
		&=(\varepsilon(x), d(b))(h(a^{-1}), d(a^2)).
	\end{align*}
	Thus ${\rm L}(A)$ is a central subgroup of $\GC(A)$. Now set 
	\[
	\PGC(A):=\GC(A)/{\rm L}(A).
	\]
	Then we have the extension
	\[
	1 \arr \U(A) \arr \PGC(A) \arr \PGE_2(A) \arr 1.
	\]

	\begin{thm}\label{E1S}
		For any commutative ring $A$, we have the exact sequence
		\[
		H_2(\PGE_2(A),\z) \arr \bigg(\displaystyle{\frac{\K_2(2, A)}{\CC(2, A)}}\bigg)_{\PGE_2(A)}^\ab
		\!\!\!\!\!\!\overset{\kappa}{\arr} A_\aa  \arr H_1(\PGE_2(A),\z) \arr \GG_A \arr 1,
		\]
		where the map on the right has a splitting $\GG_A \arr H_1(\PGE_2(A),\z)$.
	\end{thm}
	\begin{proof}
		From the Lyndon/Hochschild-Serre spectral sequence associated to the above  extension
		we obtain the five term exact sequence
		\[
		H_2(\PGC(A),\z) \arr H_2(\PGE_2(A),\z) \arr H_1(\U(A),\z)_{\PGE_2(A)} \arr H_1(\PGC(A),\z) 
		\]
		\[
		\arr H_1(\PGE_2(A),\z) \arr 0
		\]
		(see \cite[Corollary 6.4, Chp. VII]{brown1994}). By Theorem \ref{hut}, 
		\[
		H_1(\U(A),\z)_{\PGE_2(A)}\simeq \displaystyle{\bigg(\frac{\K_2(2, A)}{\CC(2, A)}}\bigg)_{\PGE_2(A)}^\ab.
		\]
		We prove that 
		\[
		H_1(\PGC(A),\z)\simeq \GG_A\oplus A_\aa.
		\]
		From the split extension $1 \arr \CC(A) \arr \GC(A) \arr \aa \arr 1$ we get the split exact sequence
		\[
		0 \arr H_1(\CC(A),\z)_\aa \arr H_1(\GC(A),\z) \arr \aa \arr 1.
		\]
		We know that $H_1(\CC(A),\z)\simeq A/M$ (see the proof of \cite[Theorem 3.1]{B-E2024}), where $M$ is 
		the additive subgroup of $A$ generated by $x(a^2-1)$ and $3(b+1)(c+1)$, $x\in A$, $a,b,c \in \aa$. This 
		map is induced by the map 
		\[
		\CC(A) \arr A/M, \ \ \ \ \epsilon(a)\mapsto a-3,
		\]
		which takes $h(a)$ to $-3(a-1)$. Following the action of $\aa$ on $\CC(A)$, we see that the action of 
		$\aa$ on $A/M$, through $ H_1(\CC(A),\z)$, is given by
		\[
		a. \overline{x}:=-\overline{(a^{-1}x+3)(a-1)}.
		\]
		Therefore
		\[
		H_1(\CC(A),\z)_\aa\simeq (A/M)_\aa \simeq A/\{ y(a-1):y\in A, a\in \aa\}=A_\aa.
		\]
		Now it follows from the above exact sequence that
		\[
		H_1(\GC(A),\z)\simeq \aa \oplus A_\aa.
		\]
		From the extension $1\arr {\rm L}(A) \arr \GC(A) \arr \PGC(A) \arr 1$ we obtain the exact sequence
		\[
		H_1({\rm L}(A),\z) \arr H_1(\GC(A),\z) \arr H_1(\PGC(A),\z) \arr 0.
		\]
		Now under the isomorphism $H_1(\GC(A),\z)\simeq \aa \oplus A_\aa$, we have
		\[
		(h(a^{-1}),d(a^2)) \mapsto (a^2, -3(a-1))=(a^2, 0).
		\]
		Thus
		\[
		H_1(\PGC(A),\z)\simeq \GG_A \oplus A_\aa.
		\]
		From the above arguments one sees that the above isomorphism is induced by the map
		\[
		\PGC(A)\arr \GG_A \oplus A_\aa, \ \ \ \ \ (\varepsilon(x), d(a)) \mapsto (\lan a\ran, \overline{x-1}).
		\]
		Composing this with $\PGSt(2, A) \arr \PGC(A)$, we get the map
		\[
		\alpha: \PGSt(2, A) \arr \GG_A \oplus A_\aa, \ \ \ \ \ (x_{ij}(t), d(a)) \mapsto (\lan a\ran, \overline{t}).
		\]
		It follows from this that the restriction of $\alpha$ to $\St(2,A)$ is given by
		\[
		\alpha: \St(2, A) \arr \GG_A \oplus A_\aa, \ \ \ \ \ 
		x_{ij}(t) \mapsto (\lan 1\ran, \overline{t})=(\lan 1\ran, f(x_{ij}(t))).
		\]
		This shows that the map
		\[
		\bigg(\displaystyle{\frac{\K_2(2, A)}{\CC(2, A)}}\bigg)_{\PGE_2(A)}^\ab \arr \GG_A \oplus A_\aa
		\]
		is given by $x \mapsto (\lan 1 \ran, \kappa(x))$. The determinant $\det: \PGE_2(A) \arr \GG_A$ induces the map 
		$\det_\ast: H_1(\PGE_2(A),\z) \arr \GG_A$. This map splits the composition 
		\[
		\GG_A\arr \GG_A\oplus A_\aa \arr H_1(\PGE_2(A),\z).
		\]
		All these give the exact sequence of the theorem.
	\end{proof} 
	
	Let $a, b\in A$ be any two elements such that $1-ab\in \aa$. We define
	\[
	\lan a,b\ran_{ij}:=x_{ji}\bigg(\frac{-b}{1-ab}\bigg)x_{ij}(-a)x_{ji}(b)x_{ij}\bigg(\frac{a}{1-ab}\bigg)h_{ij}(1-ab)^{-1}.
	\]
	It is easy to verify that  $\lan a,b\ran_{ij}\in \K_2(2, A)$. This element is called a {\it Dennis-Stein symbol}.
	If $u,v\in\aa$, then 
	\[
	\displaystyle\{u,v\}_{ij}=\bigg\lan u, \frac{1-v}{u}\bigg\ran_{ij}=\bigg\lan \frac{1-u}{v},v\bigg\ran_{ij}.
	\]
	Hence Dennis-Stein symbols generalize Steinberg symbols.
	
	\begin{cor}\label{DS}
		If $\K_2(2,A)$ is generated by Dennis-Stein symbols, then 
		\[
		H_1(\PGE_2(A),\z) \simeq \GG_A \oplus A_\aa.
		\]
	\end{cor}
	\begin{proof}
		It is easy to check that $\kappa(\lan a,b\ran_{ij})=0$. This implies that $\kappa=0$. Now the claim follows from the
		above theorem.
	\end{proof}
	
	\begin{cor}
		If $2\in \aa$, then $H_1(\PGE_2(A), \z)\simeq \GG_A$.
	\end{cor}
	\begin{proof}
		Since $2\in \aa$, $1=2-1\in \lan a-1:a\in \aa\ran$. Thus $A_\aa=0$ and the claim follows from the above theorem.
	\end{proof}
	
	\begin{exa}\label{exa-H1}
		(i) If $A$ is local with maximal ideal $\mmm_A$, then $A_\aa=0$ when $|A/\mmm_A|\neq 2$ and 
		$A_\aa=A/\mmm_A\simeq \F_2$ when $|A/\mmm_A|= 2$. Thus
		\[
		H_1(\PGL_2(A),\z)\simeq 
		\begin{cases}
			\GG_A & \text{if $|A/\mmm_A|\neq 2$} \\
			\GG_A \oplus \z/2  & \text{if $|A/\mmm_A|= 2$}.
		\end{cases}
		\]
		\par (ii) Let $A$ be a semilocal ring such that 
		none of $\z/2 \times \z/2$ and $\z/6$ is a direct factor of $A/J(A)$. Then $A$ is an universal $\GE_2$-ring
		and so by the above theorem
		\[
		H_1(\PGL_2(A),\z)\simeq \GG_A\oplus A_\aa.
		\]
		\par (iii) The ring of integers $\z$ is an universal $\GE_2$-ring \cite[Example 6.12]{hutchinson2022}. Since  
		$\z_{\z^\times}\simeq \z/2$ by the above theorem we have 
		\[
		H_1(\PGL_2(\z),\z)\simeq \z/2 \oplus \z/2.
		\]
		\par (iv) Let $m$ be a square free integer. The ring $A_m:=\z[\frac{1}{m}]$ is a $\GE_2$-ring and 
		\[
		(A_m)_{A_m^\times}\simeq \begin{cases}
			0 & \text{if $2\mid m$}\\
			\z/2 & \text{if $2\nmid m$}.
		\end{cases}
		\]
		If $m$ is odd, then the inclusion $A_m \subseteq \z_{(2)}=\{a/b: a, b \in \z, 2\nmid b\}$
		induces the commutative diagram with exact rows
		\[
		\begin{tikzcd}
			& (A_m)_{A_m^\times} \ar[r] \ar[d, "\simeq"] & H_1(\PGL_2(A_m),\z) \ar[r] \ar[d] & \GG_{A_m} \ar[r] \ar[d] & 0\\
			0 \ar[r] & (\z_{(2)})_{\z_{(2)}^\times} \ar[r] & H_1(\PGL_2((\z_{(2)}),\z) \ar[r]  & \GG_{\z_{(2)}} \ar[r] & 0.
		\end{tikzcd}
		\]
		Since the left map is an isomorphism, it follows from this diagram that the map $(A_m)_{A_m^\times} \arr H_1(\PGL_2(A_m),\z)$ 
		is injective. Therefore
		\[
		H_1(\PGL_2(A_m),\z)\simeq \begin{cases}
			\GG_{A_m} & \text{if $2\mid m$}\\
			\GG_{A_m}\oplus \z/2 & \text{if $2\nmid m$}.
		\end{cases}
		\]
		This implies that the map $\kappa$ is trivial and it follows from this and the above theorem that for any 
		$m$, the natural map
		\[
		\begin{array}{c}
			H_2(\PGL_2(\z[\frac{1}{m}]),\z) \arr 
			\bigg(\displaystyle{\frac{\K_2(2, \z[\frac{1}{m}])}{\CC(2, \z[\frac{1}{m}])}}\bigg)_{\PGL_2(\z[\frac{1}{m}])}^\ab
		\end{array}
		\]
		is surjective. The $K$-group $\K_2(2, \z[\frac{1}{m}])$ has been studied in many articles. For example 
		see \cite{morita2020} and \cite{hutchinson2016} and their references.
	\end{exa}
	
	\section{The 
		\texorpdfstring{$\GE_2$}{Lg}-unimodular vectors}\label{sec4}
	
	A column vector ${\pmb u}={\mtx {u_1} {u_2}}\in A^2$ is called {\it unimodular} if there exists 
	${\pmb v}={\mtx {v_1} {v_2}}\in A^2$ such that $({\pmb u},{\pmb v}):={\mtxx {u_1} {v_1} {u_2} {v_2}}\in\GL_2(A)$ 
	and it is called $\GE_2$-{\it unimodular} if $({\pmb u},{\pmb v}) \in \GE_2(A)$. 
	
	\begin{lem}\label{GE}
		If $\pmb{v}\in A^2$ is $\GE_2$-unimodular and if $M = (\pmb{v}, \pmb{w}) \in \GL_2(A)$, then 
		$M \in \GE_2(A)$ and $\pmb{w}$ is $\GE_2$-unimodular.
	\end{lem}
	\begin{proof}
		By definition, $\pmb{v}$ is $\GE_2$-unimodular if there exists $N \in \GE_2(A)$ with $N\pmb{e_1} =\pmb{v}$, where 
		${\pmb e_1}:={\mtx 1 0}$. Hence $\pmb{v}$ is $\GE_2$-unimodular if and only if there exists $P \in \GE_2(A)$ with 
		$P\pmb{v} = \pmb{e_1}$. Thus
		\[
		PM = (\pmb{e_1}, P\pmb{w}).
		\]
		It is clear that a matrix of the form $(\pmb{e_1}, \pmb{u})$ is invertible if
		and only if it lies in $\GE_2(A)$.
	\end{proof}
	
	For any non-negative integer $n$, let $L_n(A^2)$ be the free abelian group generated 
	by the set of all $(n+1)$-tuples $(\lan{\pmb v_0}\ran, \dots, \lan{\pmb v_n}\ran)$ of unimodular vectors 
	${\pmb v_i}\in A^2$ such that for any $i\neq j$, the matrix $({\pmb v_i},{\pmb v_j})$ 
	is invertible. Note that for a vector $\pmb{v}\in A^2$, $\lan \pmb{v}\ran$ means the 
	equivalence class up to multiplication by a unit, i.e. $\lan \pmb{v}\ran=\pmb{v}\aa$.
	
	We consider $L_n(A^2)$ as a left $\PGL_2(A)$-module in a natural way. If necessary, we convert this 
	action to a right action by the definition $m.g:=g^{-1}m$. We define the $n$-th differential operator
	\[
	\partial_n^L : L_n(A^2) \arr L_{n-1}(A^2), \ \ n\ge 1,
	\]
	as an alternating sum of face operators which throws away the $i$-th component of generators. Then 
	\[
	L_\bullet(A^2): \  \cdots \larr  L_2(A^2) \overset{\partial_2^L}{\larr} L_1(A^2) 
	\overset{\partial_1^L}{\larr} L_0(A^2) \arr 0
	\]
	is a complex. This complex has been studied extensively in \cite{hutchinson2022} and \cite{C-H2022}.
	
	Let $Y_n(A^2)$ be the free abelian subgroup of $L_n(A^2)$ generated by the set of all $(n+1)$-tuples 
	$(\lan{\pmb v_0}\ran, \dots, \lan{\pmb v_n}\ran)$ of $\GE_2$-unimodular vectors. Thus $Y_\bullet(A^2)$
	is a $\PGE_2(A)$-subcomplex of $L_\bullet(A^2)$. We say that $Y_\bullet(A^2)$ (resp. $L_\bullet(A^2)$) is 
	{\it exact in dimension} $k$ if $H_k(Y_\bullet(A^2))=0$ (resp. $H_k(L_\bullet(A^2))=0$). 
	
	For a subgroup $H$ of a group $G$ and any $H$-module $M$, let $\Ind_H^GM:=\z[G]\otimes_H M$. This extension
	of scalars is called {\it induction} from $H$ to $G$.
	
	\begin{lem}\label{ind1}
		The natural inclusion $Y_\bullet(A^2)\arr L_\bullet(A^2)$ induces the isomorphism
		\[
		L_\bullet(A^2)\simeq \Ind_{\PGE_2(A)}^{\PGL_2(A)} Y_\bullet(A^2).
		\]
	\end{lem}
	\begin{proof}
		Clearly 
		\[
		\phi_\bullet: \z[\PGL_2(A)] \otimes_{\PGE_2(A)} Y_\bullet(A^2) \arr L_\bullet(A^2)
		\]
		given by
		\[
		g\otimes (\lan\pmb{v_0}\ran, \lan\pmb{v_1}\ran, \dots, \lan\pmb{v_n}\ran) \mapsto 
		(\lan g\pmb{v_0}\ran, \lan g\pmb{v_1}\ran, \dots, \lan g\pmb{v_n}\ran),
		\]
		is a well-defined morphism of complexes of $\z[\PGL_2(A)]$-modules. In fact $\phi_\bullet$ is an 
		isomorphism with the inverse morphism
		\[
		\psi_\bullet: L_\bullet(A^2) \arr \z[\PGL_2(A)] \otimes_{\PGE_2(A)} Y_\bullet(A^2)
		\]
		defined by
		\[
		(\lan\pmb{v_0}\ran, \lan\pmb{v_1}\ran, \dots, \lan\pmb{v_n}\ran) \mapsto g \otimes 
		(\lan\pmb{e_1}\ran, \lan g^{-1}\pmb{v_1}\ran, \dots, \lan g^{-1}\pmb{v_n}\ran),
		\]
		where $g\pmb{e_1}=\pmb{v_0}$. Note that by Lemma \ref{GE},  $\psi_\bullet$ is well-defined.
	\end{proof}
	
	\begin{cor}\label{ind2}
		For any non-negative integer $n$, 
		\[
		H_n(L_\bullet(A^2))\simeq \Ind_{\PGE_2(A)}^{\PGL_2(A)} H_n(Y_\bullet(A^2)).
		\]
	\end{cor}
	\begin{proof}
		Since the functor $\Ind_{\PGE_2(A)}^{\PGL_2(A)}$ is exact on the category of $\z[\PGE_2(A)]$-modules (since
		$\z[\PGL_2(A)]$ is a free $\z[\PGE_2(A)]$-module), the claim follows from the previous lemma.
	\end{proof}
	
	The group $\SL_2(A)$ (resp. $\Ee_2(A)$) acts transitively on the sets of generators of $L_0(A^2)$ 
	(resp. $Y_0(A^2)$). Let
	\[
	{\pmb\infty}:=\lan {\pmb e_1}\ran, \ \ \  {\pmb 0}:=\lan {\pmb e_2}\ran , \ \ \  
	{\pmb a}:=\lan {\pmb e_1}+ a{\pmb e_2}\ran, \ \ \ a\in \aa,
	\]
	where ${\pmb e_1}:={\mtx 1 0}$ and $ {\pmb e_2}:={\mtx 0 1}$. Let $\epsilon: L_0(A^2) \arr \z$ be 
	defined by $\sum_i n_i(\lan {\pmb v_{0,i}}\ran) \mt \sum_i n_i$. We denote the restriction  
	$\epsilon|_{Y_0(A^2)}:Y_0(A^2) \arr \z$ again by $\epsilon$.
	
	\begin{prp}[Hutchinson]\label{GE2C-1}
		For any commutative  ring $A$, $H_0(Y_\bullet(A^2))\overset{\epsilon}{\simeq} \z$. In other words, the complex
		\[
		Y_1(A^2) \overset{\partial_1^Y}{\larr} Y_0(A^2)\overset{\epsilon}{\arr}\z\arr 0
		\]
		always is exact. Moreover, $H_0(L_\bullet(A^2))\simeq \Ind_{\PGE_2(A)}^{\PGL_2(A)} \z$.
	\end{prp} 
	\begin{proof}
		Here we follow the proof of \cite[Theorem 3.3]{hutchinson2022}.
		Clearly $\epsilon:Y_0(A^2)\arr \z$ is surjective. Let $X\in \ker(\epsilon)$. We may 
		assume $X=(\lan \pmb{u}\ran )-(\lan \pmb{v}\ran )$. Since $\GE_2(A)$ acts transitively 
		on the generators of $Y_0(A^2)$, we may assume $X=\pmb{\infty}-E\pmb{\infty}$, where 
		$E\in \GE_2(A)$. For any $x\in A$ and $a, b\in \aa$, we have 
		\[
		E(x)\diag(a, b)=\diag(b, a)E(b^{-1}xa).
		\]
		Thus any element of $\GE_2(A)$ can be written as product $E'D'$, where $D'\in {\rm D}_2(A)$ and 
		$E'\in \Ee_2(A)$. Since $D'\pmb{\infty}=\pmb{\infty}$, we may assume that $E\in \Ee_2(A)$. Let
		\[
		E=E(a_1)^{c_1}\cdots E(a_n)^{c_n},
		\]
		where $c_i \in \{1,-1\}$.
		If $E_i:=E(a_1)^{c_1}\cdots E(a_i)^{c_i}$ for $1\leq i \leq n$, and $E_0=I_2$, then 
		\[
		Y:=\sum_{i=1}^{n} (E_{i}\pmb{\infty}, E_{i-1} \pmb{\infty})\in Y_1(A^2)
		\]
		and $\epsilon(Y)=X$. This proves our claim.
	\end{proof}
	
	Over the class of local rings we have the following result of Hutchinson.
	
	\begin{prp}[Hutchinson]\label{<k}
		Let $A$ be a local ring. Then the complex $Y_\bullet(A^2) \overset{\epsilon}{\arr} \z$ is exact in 
		dimension $<|A/\mmm_A|$.
	\end{prp}
	\begin{proof}
		See \cite[Lemma 3.21]{hutchinson2017}.
	\end{proof}
	
	It follows from Lemma \ref{ind1} and Shapiro’s Lemma that the inclusion $Y_\bullet(A^2)\arr L_\bullet(A^2)$
	induces isomorphisms of the homology groups
	\[
	H_q(\PGE_2(A), Y_p(A^2)) \simeq H_q(\PGL_2(A), L_p(A^2))
	\]
	(for all $ p,q$), which occur in the spectral sequences in Section \ref{sec6} and in \cite[\S2]{C-H2022},
	\cite[\S7]{hutchinson2022}. For any $n\geq 0$, let
	\[
	Z_n^{\GE_2}(A^2):=\ker(\partial_n^Y), \ \ \ \ \ Z_n^{\GL_2}(A^2):=\ker(\partial_n^L).
	\]
	Now from the above isomorphism and the fact that the functor $\Ind_{\PGE_2(A)}^{\PGL_2(A)}$ is exact on the 
	category of $\z[\PGE_2(A)]$-modules, we obtain the isomorphism
	\[
	Z_n^{\GL_2}(A^2) \simeq \Ind_{\PGE_2(A)}^{\PGL_2(A)} Z_n^{\GE_2}(A^2).
	\]
	From this we obtain the following lemma.
	
	\begin{lem}\label{ind3}
		For all $p,q\geq 0$, we have the isomorphism 
		\[
		H_q(\PGE_2(A),  Z_p^{\GE_2}(A^2)) \simeq H_q(\PGL_2(A), Z_p^{\GL_2}(A^2)).
		\]
	\end{lem}
	
	\section{ The homology group \texorpdfstring{$H_1(Y_\bullet(A^2))$}{Lg}}\label{sec5}
	
	Let $\Gamma(A)$ be the graph of unimodular rows introduced and studied in \cite[\S2]{hutchinson2022} and let 
	$\Gamma^{\GE}(A)$ be the analogous graph of $\GE_2$-unimodular rows. Then Lemma \ref{GE} shows that 
	$\Gamma^{\GE}(A)$ is precisely the path component of $\pmb{\infty} = \lan \pmb{e_1} \ran$ in $\Gamma(A)$. 
	Furthermore, the transitive action of $\PGL_2(A)$ on $\Gamma(A)$ shows that it decomposes into 
	homeomorphic path components
	\[
	\Gamma(A)= \bigsqcup_{g\in \PGL_2(A)/\PGE_2(A)} g.\Gamma^{\GE}(A).
	\]
	If we now let $Y(A)$ denote the clique complex of $\Gamma(A)$ as in \cite[\S2]{hutchinson2022} and if we let 
	$Y^{\GE}(A)$ denote the clique complex of $\Gamma^{\GE}(A)$, then it follows that
	\[
	Y(A)= \bigsqcup_{g\in \PGL_2(A)/\PGE_2(A)} g.Y^{\GE}(A).
	\]
	Taking geometric realizations it again follows that
	\[
	|Y(A)|= \bigsqcup_{g\in \PGL_2(A)/\PGE_2(A)} g.|Y^{\GE}(A)|
	\]
	and that $|Y^\GE(A)|$ is the path component at $\pmb{\infty}$ of $|Y(A)|$. In particular, it follows that the 
	inclusion $|Y^\GE(A)|\arr |Y(A)|$ induces the following result:
	
	\begin{prp}[Hutchinson]\label{Pi1}
		For any commutative ring $A$, we have the isomorphism
		\[
		\pi_1(|Y^\GE(A)|, \pmb{\infty})= \pi_1(|Y(A)|, \pmb{\infty})\simeq \frac{\K_2(2,A)}{\CC(2, A)}.
		\]
	\end{prp} 
	\begin{proof}
		See \cite[Theorem 6.9]{hutchinson2022}.
	\end{proof}
	Since the space $|Y^\GE(A)|$ is path-connected, it follows from the above theorem that
	\[
	H_1(|Y^\GE(A)|, \z)\simeq \pi_1(|Y^\GE(A)|, \pmb{\infty})\simeq \bigg(\frac{\K_2(2,A)}{\CC(2, A)}\bigg)^\ab.
	\]
	Let $\Delta^\GE(A)$ denote the standard ordered chain complex of the simplicial complex $Y^\GE(A)$. As in 
	\cite[\S7]{hutchinson2022}, the complex $Y_\bullet(A^2)$ in the current article is the complex of non-degenerate ordered 
	simplices of the simplicial complex $Y^\GE(A)$ and the natural map of complexes $Y_\bullet(A^2) \arr \Delta^\GE(A)$ 
	induces an isomorphism on first homology groups. Thus
	\[
	H_1(Y_\bullet(A^2)) \simeq H_1(\Delta^\GE(A),\z) = H_1(|Y^\GE(A)|,\z) \simeq \bigg(\frac{K_2(2, A)}{\CC(2, A)}\bigg)^\ab.
	\]
	Thus we have:
	
	\begin{thm}[Hutchinson]\label{GE2CC-1}
		For any commutative ring $A$, we have the isomorphisms
		\[
		H_1(Y_\bullet(A^2))\simeq \bigg(\displaystyle{\frac{\K_2(2, A)}{\CC(2, A)}}\bigg)^\ab, \ \ \ \ \
		H_1(L_\bullet(A^2))\simeq \Ind_{\PGE_2(A)}^{\PGL_2(A)}\bigg(\displaystyle{\frac{\K_2(2, A)}{\CC(2, A)}}\bigg)^\ab.
		\]
		In particular if $A$ is universal for $\GE_2$, then $Y_\bullet(A^2) \overset{\epsilon}{\arr} \z$ and 
		$L_\bullet(A^2) \overset{\epsilon}{\arr} \z$ are exact in  dimension $1$.
	\end{thm}
	
	\begin{rem}
		In  \cite[Theorem 7.2]{hutchinson2022}, Hutchinson states that 
		$H_1(L_\bullet(A^2))\simeq \bigg(\displaystyle{\frac{\K_2(2, A)}{\CC(2, A)}}\bigg)^\ab$.
		In fact this is only valid when the space $Y(A)$ is path-connected; i.e., when $A$ is a 
		$\GE_2$-ring. Theorem \ref{GE2CC-1} above gives a corrected statement valid for all 
		commutative rings.
	\end{rem}
	
	\section{The main spectral sequence}\label{sec6}
	
	Let $A$ be a commutative ring. The group $\PGE_2(A)$ acts naturally on $Z_i^{\GE_2}(A^2)$. 
	By Proposition~\ref{GE2C-1}, the sequence of $\PGE_2(A)$-modules
	\begin{equation*}\label{comp1}
		0 \arr Z_1^{\GE_2}(A^2) \overset{\inc}{\arr} Y_1(A^2) \overset{\partial_1^Y}{\arr} Y_0(A^2) \arr \z \arr 0
	\end{equation*}
	is exact. Let $E_\bullet(A^2)$ be the sequence
	\begin{equation*}
		0 \arr Z_1^{\GE_2}(A^2) \overset{\inc}{\arr} Y_1(A^2) \overset{\partial_1^Y}{\arr} Y_0(A^2) \arr 0
	\end{equation*}
	and $B_\bullet(\PGE_2(A))\arr \z$ be the bar resolution of $\PGE_2(A)$ over $\z$ \cite[Chap.I, \S 5]{brown1994}. 
	Let $D_{\bullet,\bullet}$ be the double complex $B_\bullet(\PGE_2(A))\otimes_{\PGE_2(A)} E_\bullet(A^2)$.
	From this double complex we obtain the first quadrant spectral sequence 
	\[
	E^1_{p,q}=\left\{\begin{array}{ll}
		H_q(\PGE_2(A),Y_p(A^2)) & p=0,1\\
		H_q(\PGE_2(A),Z_1^{\GE_2}(A^2)) & p=2\\
		0 & p>2
	\end{array}
	\right.
	\Longrightarrow H_{p+q}(\PGE_2(A),\z)
	\]
	(see \cite[\S5, Chap. VII]{brown1994}). 
	
	The group $\PGE_2(A)$ acts transitively on the set of generators of $Y_i(A^2)$ for $i=0,1$. 
	We choose $({\pmb \infty})$ and $({\pmb \infty} ,{\pmb 0})$ as representatives of the orbit of 
	the generators of $Y_0(A^2)$ and $Y_1(A^2)$, respectively. Then
	\begin{equation}\label{shapiro}
		Y_0(A^2)\simeq \Ind _{\PB_2(A)}^{\PGE_2(A)}\z, \ \ \ \ \   Y_1(A^2)\simeq \Ind _{\PT_2(A)}^{\PGE_2(A)}\z,
	\end{equation}
	where 
	\[
	\PB_2(A):=\stabe_{\PGE_2(A)}({\pmb \infty})=\Bigg\{\begin{pmatrix}
		a & b\\
		0 & d
	\end{pmatrix}:a,d\in \aa, b\in A\bigg\}/\aa I_2,
	\]
	\[
	\PT_2(A):=\stabe_{\PGE_2(A)}({\pmb \infty},{\pmb 0})=\Bigg\{\begin{pmatrix}
		a & 0\\
		0 & d
	\end{pmatrix}:a, d\in \aa\bigg\}/\aa I_2.
	\]
	Note that $\PT_2(A)\overset{\simeq}{\larr} \aa$, which is given by $\diag(a,d)\mt ad^{-1}$. The inverse
	of this map is 
	\[
	\aa \overset{\simeq}{\larr} \PT_2(A), \ \ \ \ \ a\mt \diag(a,1)=\diag(1,a^{-1}). 
	\]
	Usually in our calculations we identify $\PT_2(A)$ with $\aa$. The group $\PB_2(A)$ sits in the split extension
	\[
	1 \arr N_2(A) \arr \PB_2(A) \arr \PT_2(A) \arr 1,
	\]
	where 
	\[
	N_2(A):=\Bigg\{\begin{pmatrix}
		1 & b\\
		0 & 1
	\end{pmatrix}: b\in A\bigg\}\simeq A.
	\]
	So we have the split extension $0 \arr A \arr \PB_2(A) \arr \aa \arr 1$,
	where the action of $\aa$ on $A$ is given by $a.x:=ax$. This implies that
	\[
	H_0(\aa,A)=A_\aa=A/\lan a-1: a\in \aa\ran.
	\]
	
	By Shapiro's lemma, applied to (\ref{shapiro}), we have
	\[
	E_{0,q}^1 \simeq H_q(\PB_2(A),\z), \  \  \ 
	E_{1,q}^1 \simeq H_q(\PT_2(A),\z).
	\]
	In particular, $E_{0,0}^1\simeq \z\simeq E_{1,0}^1$. Moreover
	\[
	d_{1, q}^1=\sigma_\ast - \inc_\ast,
	\]
	where $\sigma: \PT_2(A) \arr \PB_2(A)$ is given by $\sigma(X)= wX w^{-1}$ for 
	$w={\mtxx 0 1 {-1} 0}$. This easily implies that  $d_{1,0}^1$ is trivial, $d_{1,1}^1$ is 
	induced by the map $\PT_2(A) \arr \PB_2(A)$, $X\mt X^{-2}$, and 
	\[
	d_{1,2}^1: H_2(\PT_2(A),\z) \arr H_2(\PB_2(A),\z)
	\]
	is trivial. 
	
	\section{The scissors congruence group}\label{sec7}
	
	Following Coronado and Hutchinson \cite[\S3]{C-H2022} we define the {\it scissors congruence group} 
	of $A$ as follows:
	\[
	\PP(A):=H_0(\PGL_2(A), Z_2^{\GL_2}(A^2)).
	\]
	It follows from Lemma \ref{ind3}, that
	\[
	\PP(A)\simeq H_0(\PGE_2(A), Z_2^{\GE_2}(A^2)).
	\]
	
	\begin{rem}\label{P_E=P}
		Let $A$ satisfy the condition that the complex $Y_\bullet(A^2)\arr \z$ is exact in dimension $<4$
		(for example see Proposition \ref{<k}). 
		Then $\PP(A)$ is isomorphic with the classical scissors congruence group defined in  
		Section~\ref{sec1}. In fact, from the exact sequence $Y_4(A^2) \arr Y_3(A^2) \arr Z_2^{\GE_2}(A^2) \arr 0$ we obtain the 
		exact sequence 
		\[
		Y_4(A^2)_{\PGE_2(A)} \arr Y_3(A^2)_{\PGE_2(A)} \arr \PP(A) \arr 0. 
		\]
		The orbits of the action of  $\PGE_2(A)$ on $Y_3(A)$  and $Y_4(A)$ are represented by
		\[
		[x]':=({\pmb\infty}, {\pmb 0},{\pmb 1}, \pmb{x}),
		\ \  \text{and}  \ \   
		[x,y]':= ({\pmb\infty}, {\pmb 0},{\pmb 1}, \pmb{x}, \pmb{y}),
		\ \ \ x,y,x/y\in \WW_A,
		\]
		respectively. Thus $Y_3(A^2)_{\PGE_2(A)}$ is the free abelian group generated by the symbols $[x]'$, $x\in \WW_A$ 
		and  $Y_4(A^2)_{\PGE_2(A)}$ is the free abelian group generated by the symbols
		$[x,y]'$, $x,y,x/y\in \WW_A$. It is straightforward to check that 
		\[
		\overline{\partial_4^Y}([x,y]')=
		[x]'-[y]'+\bigg[\frac{y}{x}\bigg]'-\Bigg[\frac{1-x^{-1}}{1-y^{-1}}\Bigg]'+\Bigg[\frac{1-x}{1-y}\Bigg]'.
		\]
		This proves our claim.
	\end{rem}

	\begin{lem}\label{PE}
		If $A$ satisfies the condition that $Y_\bullet(A^2)$ is exact in dimension one, then 
		\[
		\PP(A)\simeq H_1(\GE_2(A), Z_1^{\GE_2}(A^2)). 
		\]
	\end{lem}
	\begin{proof}
		Since $Y_\bullet(A^2)$ is exact in dimension one, the sequence 
		\[
		0  \arr Z_2^{\GE_2}(A^2) \arr Y_2(A^2) \arr Z_1^{\GE_2}(A^2)\arr 0
		\]
		is exact. From this we obtain the long 
		exact sequence 
		\[
		\cdots \arr H_1(\PGE_2(A), Y_2(A^2))\arr H_1(\PGE_2(A), Z_1^{\GE_2}(A^2))   \arr H_0(\PGE_2(A), Z_2^{\GE_2}(A^2))
		\]
		\[
		\arr H_0(\PGE_2(A), Y_2(A^2)) \arr H_0(\PGE_2(A), Z_1^{\GE_2}(A^2)) \arr 0.
		\]
		The group $\PGE_2(A)$ acts transitively on the generators of $Y_2(A^2)$. We choose 
		$({\pmb \infty}, {\pmb 0}, {\pmb 1})$ as representative of the orbit of the generators 
		of $Y_2(A^2)$. Then 
		\[
		Y_2(A^2)\simeq  \Ind _{\{1\}}^{\PGE_2(A)}\z. 
		\]
		Thus by Shapiro's lemma,
		\[
		H_q(\PGE_2(A),Y_2(A^2))\simeq H_q(\{1\}, \z)\simeq
		\begin{cases}
			\z &\text{if $q=0$}\\
			0 &\text{if $q> 0$}.
		\end{cases}
		\]
		The map $\partial_2^Y: Y_2(A^2) \arr Y_1(A^2)$ induces the identity map
		\[
		\overline{\partial_2^Y}: \z=Y_2(A^2)_{\PGE_2(A)} \arr Y_1(A^2)_{\PGE_2(A)}=\z.
		\]
		Since $\overline{\partial_2^Y}$ factors through $Z_1^{\GE_2}(A^2)_{\PGE_2(A)}$ and 
		$Y_2(A^2)_{\PGE_2(A)}\arr Z_1^{\GE_2}(A^2)_{\PGE_2(A)}$ is surjective, we conclude that the map 
		\[
		Y_2(A^2)_{\PGE_2(A)}\arr Z_1^{\GE_2}(A^2)_{\PGE_2(A)}
		\]
		must be an isomorphism. Therefore
		\[
		H_1(\PGE_2(A), Z_1^{\GE_2}(A^2)) \simeq H_0(\PGE_2(A), Z_2^{\GE_2}(A^2))=\PP(A).
		\]
	\end{proof}
	
	\begin{lem}\label{E21}
		The differential $d_{2,1}^1$ is trivial. In particular, if $Y_\bullet(A)$ is exact in dimension one, 
		then $E_{2,1}^2\simeq \PP(A)$.
	\end{lem}
	\begin{proof}
		Let $D_{\bullet,\bullet}'$ be the double complex $ F_\bullet \otimes_{\GE_2(A)} E_\bullet(A^2)$,
		where $F_\bullet \arr \z$ is a projective resolution of $\GE_2(A)$ over $\z$.
		From $D_{\bullet,\bullet}'$ we obtain the first quadrant spectral sequence
		\[
		{E'}^1_{p.q}=\left\{\begin{array}{ll}
			H_q(\GE_2(A),Y_p(A^2)) & p=0,1\\
			H_q(\GE_2(A),Z_1^{\GE_2}(A^2)) & p=2\\
			0 & p>2
		\end{array}
		\right.
		\Longrightarrow H_{p+q}(\GE_2(A),\z).
		\]
		The natural map $p:\GE_2(A) \arr \PGE_2(A)$ induces the morphism of spectral sequences
		\begin{equation}\label{morphism}
			\begin{tikzcd}
				{E'}_{p,q}^1\ar[r,Rightarrow] \ar[d, "p_\ast"]& H_{p+q}(\GE_2(A),\z)\ar[d, "p_\ast"]\\
				E_{p,q}^1\ar[r,Rightarrow]& H_{p+q}(\PGE_2(A),\z).
			\end{tikzcd}
		\end{equation}
		As in case of the spectral sequence $E_{\bullet, \bullet}^1$ discussed in the previous section,
		we can show that
		\[
		{E'}_{0,q}^1=H_q({\rm B}_2(A),\z), \ \ \ \ \ {E'}_{1,q}^1=H_q({\rm T}_2(A),\z),
		\]
		where
		\[
		{\rm B}_2(A):=\Bigg\{\begin{pmatrix}
			a & b\\
			0 & d
		\end{pmatrix}:a,d\in \aa, b\in A\bigg\},
		\ \ \ 
		{\rm T}_2(A):=\Bigg\{\begin{pmatrix}
			a & 0\\
			0 & d
		\end{pmatrix}:a, d\in \aa\bigg\}.
		\]
		From the above morphism of spectral sequences we have the commutative diagram
		\[
		\begin{tikzcd}
			{E'}_{2,1}^1 \ar[r, "{d'}_{2,1}^1"] \ar[d, "p_\ast"]& H_1({\rm T}_2(A),\z) \ar[r, "{d'}_{1,1}^1"]\ar[d, "p_\ast"]& 
			H_1({\rm B}_2(A),\z)\ar[d, "p_\ast"]\\
			E_{2,1}^1 \ar[r, "{d}_{2,1}^1"] & H_1(\PT_2(A),\z) \ar[r, "{d}_{1,1}^1"] & H_1(\PB_2(A),\z).
		\end{tikzcd}
		\]
		The differential ${d'}_{1,1}^1$ can be calculated similar to ${d}_{1,1}^1$. In fact 
		${d'}_{1,1}^1={\sigma'}_\ast-\inc_\ast$, where $\sigma': {\rm T}_2(A) \arr {\rm B}_2(A)$, 
		$\diag(a,b)\mapsto \diag(b,a)$. It is easy to see that
		\[
		{d'}_{1,1}^1(\diag(a,b))=\diag(ba^{-1}, ab^{-1}).
		\]
		It follows from this that $\ker({d'}_{1,1}^1)=\aa I_2$ and thus $p_\ast\circ {d'}_{2,1}^1=0$. Since the vertical maps
		are surjective \cite[Corollary 6.4, Chap. VII]{brown1994}, the differential $d_{2,1}^1$ must be trivial.
		The second part follows from the first part and Lemma \ref{PE}.
	\end{proof}
	
	\begin{thm}\label{ES}
		Let $A$ be a commutative ring which satisfies the condition that $Y_\bullet(A^2)$ is exact in dimension one. 
		Then $H_1(\PGE_2(A),\z)\simeq \GG_A\oplus A_\aa$ and we have the exact sequence
		\[
		H_3(\PGE_2(A),\z) \arr \PP(A) \arr H_2(\PB_2(A),\z) \arr H_2(\PGE_2(A),\z) \arr \mu_2(A) \arr 1.
		\]
	\end{thm}
	\begin{proof}
		Consider the composite $Y_2(A^2) \overset{\partial_2^Y}{\larr} Z_1^{\GE_2}(A^2) \overset{\inc}{\larr} Y_1(A^2)$. 
		Since $H_1(Y_\bullet(A^2))=0$, the left map is surjective. As we discussed in the proof of Lemma \ref{PE}, the map
		\[
		\overline{\partial_2^Y}:\z\simeq Y_2(A^2)_{\PGE_2(A)}  \arr Y_1(A^2)_{\PGE_2(A)}\simeq \z
		\]
		is an isomorphism. This implies that the differential 
		\[
		d_{2,0}^1=\overline{\inc}: Z_1^{\GE_2}(A^2)_{\PGE_2(A)} \arr Y_1(A^2)_{\PGE_2(A)}=\z
		\]
		is surjective. On the other hand $Y_2(A^2)_{\PGE_2(A)} \overset{\partial_2^Y}{\larr} Z_1^{\GE_2}(A^2)_{\PGE_2(A)}$ 
		is surjective. Thus $d_{2,0}^1$ is injective too and therefore 
		\[
		E_{2,0}^2=0.
		\]
		On the other hand, we have
		\[
		E_{1,1}^2\simeq \mu_2(A), \ \ \ E_{0,2}^2=H_2(\PB_2(A),\z).
		\]
		From the split extension $0\arr A \arr \PB_2(A) \arr \aa\arr 1$, we obtain the five term exact sequence
		\[
		H_2(\PB_2(A),\z) \arr H_2(\aa,\z) \arr  A_\aa \arr H_1(\PB_2(A),\z) \arr \aa \arr 1,
		\]
		(see \cite[Corollary 6.4, Chap. VII]{brown1994}). Since the above extension splits, the left side map in the above 
		exact sequence is surjective. Thus 
		\[
		H_1(\PB_2(A),\z)\simeq \aa\oplus A_\aa.
		\]
		Since under the differential 
		\[
		d_{1,1}^1:\aa\simeq H_1(\PT_2(A),\z) \arr H_1(\PB_2(A),\z)\simeq \aa\oplus A_\aa,
		\]
		the element $a\in \aa$ maps to $(a^{-2},0)$, we have $E_{0,1}^2\simeq\GG_A\oplus A_\aa$.
		Now the theorem follows from an easy analysis of the main spectral sequence. 
	\end{proof}
	
	\begin{rem}
		Let $A$ be a commutative ring. From the diagram with  exact rows
		\[
		\begin{tikzcd}
			Y_2(A^2) \ar[r,"\partial_2^Y"] \ar[d,"\partial_2^Y"] & Z_1^{\GE_2}(A^2) \ar[r]\ar[d, "\inc"] & 
			H_1(Y_\bullet(A^2)) \ar[r] & 0\\
			Y_1(A^2) \ar[r, equal ] & Y_1(A^2) & &
		\end{tikzcd}
		\]
		we obtain the commutative diagram with exact rows
		\[
		\begin{tikzcd}
			\z= Y_2(A^2)_{\PGE_2(A)} \ar[r,"\overline{\partial_2^Y}"] \ar[d,"\id_\z=\overline{\partial_2^Y}"] & 
			Z_1^{\GE_2}(A^2)_{\PGE_2(A)} \ar[r]\ar[d, "d_{2,0}^1=\overline{\inc}"] & H_1(Y_\bullet(A^2))_{\PGE_2(A)}\ar[r] & 0\\
			\z \ar[r, equal ] & \z & &
		\end{tikzcd}
		\]
		(see the proof of the above theorem). Thus by the Snake lemma we have 
		\[
		E_{2,0}^2\simeq H_1(Y_\bullet(A^2))_{\PGE_2(A)}.
		\]
		Since
		\[
		E_{1,1}^2\simeq \mu_2(A), \ \ \ E_{0,2}^2=H_2(\PB_2(A),\z), \ \ \ E_{0,1}^2\simeq\GG_A\oplus A_\aa,
		\]
		by an easy analysis of the  main spectral sequence, we obtain the exact sequence
		\[
		H_2(\PGE_2(A),\z) \arr H_1(Y_\bullet(A^2))_{\PGE_2(A)} \arr \GG_A\oplus A_\aa  \arr H_1(\PGE_2(A),\z) \arr 0.
		\]
		Combining this with Theorem \ref{GE2CC-1} we obtain the exact sequence 
		\[
		H_2(\PGE_2(A),\z) \arr \bigg(\displaystyle{\frac{\K_2(2, A)}{C(2, A)}}\bigg)_{\PGE_2(A)}^\ab
		\!\!\!\!\!\!\arr A_\aa  \arr H_1(\PGE_2(A),\z) \arr \GG_A \arr 1.
		\]
		We believe that this exact sequence coincides with the exact sequence of Theorem \ref{E1S}. It seems very difficult 
		to describe the map $\bigg(\displaystyle{\frac{\K_2(2, A)}{C(2, A)}}\bigg)_{\PGE_2(A)}^\ab \arr A_\aa$ in the above 
		exact sequence using the differentials of the spectral sequence, while it was reasonably easy to describe a similar 
		map in Theorem \ref{E1S}.
	\end{rem}
	
	\begin{cor}\label{H1-uni}
		If $A$ is universal for $\GE_2$, then 
		we have the exact sequence
		\[
		H_3(\PGE_2(A),\z) \arr \PP(A) \arr H_2(\PB_2(A),\z) \arr H_2(\PGE_2(A),\z) \arr \mu_2(A) \arr 1.
		\]
	\end{cor}
	\begin{proof}
		Since $A$ is universal for $\GE_2$, by Theorem \ref{GE2CC-1} we have 
		\[
		H_1(Y_\bullet(A^2))\simeq \bigg(\displaystyle{\frac{\K_2(2, A)}{C(2, A)}}\bigg)^\ab=0.
		\]
		Now the claim follows from the above theorem. 
	\end{proof}
	
	\begin{exa}
		Let $A$ be a semilocal ring such that 
		none of $\z/2 \times \z/2$ and $\z/6$ is a direct factor of $A/J(A)$. Then $A$ is a universal $\GE_2$-ring 
		and thus by the above corollary 
		we have the exact sequence
		\[
		H_3(\PGL_2(A),\z) \arr \PP(A) \arr H_2(\PB_2(A),\z) \arr H_2(\PGL_2(A),\z) \arr \mu_2(A) \arr 1.
		\]
	\end{exa}
	
	\section{The homology groups of \texorpdfstring{$\PB_2(A)$}{Lg}}\label{sec8}
	
	Let study the Lyndon/Hochschild-Serre spectral sequence associated to the split extension
	\[
	1 \arr N_2(A) \arr \PB_2(A) \arr \PT_2(A) \arr 1.
	\]
	This is the extension $0 \arr A \arr \PB_2(A) \arr \aa \arr 1$. Thus we have the spectral sequence
	\[
	\EE_{p,q}^2=H_p(\aa, H_q(A, \z)) \Longrightarrow H_{p+q}(\PB_2(A),\z).
	\]
	We showed in the proof of Theorem \ref{ES} that
	\[
	H_1(\PB_2(A),\z) \simeq \aa \oplus A_\aa \simeq H_1(\PT_2(A),\z) \oplus A_\aa.
	\]
	Recall that $A_\aa=H_0(\aa, A)=A/\lan a-1: a\in \aa\ran$.
	
	\begin{lem}\label{lem-sus}
		Let $G$ be an abelian group, $A$ a commutative ring, $M$ an $A$-module and 
		$\varphi: G \arr A^\times$ a homomorphism of groups which turns $A$ and $M$ 
		into $G$-modules.  If $H_0(G,A)=0$, then for any  $n\geq 0$, $H_n(G,M)=0$.
	\end{lem}
	\begin{proof}
		See \cite[Lemma~1.8]{nes-suslin1990}.
	\end{proof}
	
	\begin{cor}\label{H0A}
		If $A_\aa=0$, then $H_n(\aa,A)=0$ for any $n\geq 0$.
	\end{cor}
	\begin{proof}
		Use the above lemma by considering $\varphi=\id_\aa: \aa \arr \aa$.
	\end{proof}
	
	\begin{cor}\label{2-inv}
		If $2\in \aa$, then $H_n(\aa,A)=0$ for any $n\geq 0$.
	\end{cor}
	\begin{proof}
		If $2\in \aa$, then $A_\aa=0$. Now use the previous corollary.
	\end{proof}
	
	\begin{lem}
		If $A_\aa=0$, then $H_2(\PB_2(A),\z) \simeq H_2(\PT_2(A),\z) \oplus H_2(A,\z)_\aa$.
	\end{lem}
	\begin{proof}
		By Corollary \ref{H0A}, we have  $H_1(\aa,A)=0$. Now the claim follows from an easy analysis
		of the above spectral sequence.
	\end{proof}
	
	\begin{prp}\label{1/2}
		If $A$ is a subring of $\q$, then for any $n\geq 0$, 
		\[
		H_n(\PB_2(A),\z)\simeq H_n(\PT_2(A),\z) \oplus H_{n-1}(\aa, A).
		\]
		In particular if $2\in \aa$, then $H_n(\PT_2(A),\z) \simeq H_n(\PB_2(A),\z)$.
	\end{prp}
	\begin{proof}
		It is well known that any finitely generated subgroup of $\q$ is cyclic. Thus the additive group
		$A$ is a direct limit of infinite cyclic groups. Since $H_n(\z,\z)=0$ for any
		$n\geq  2$ \cite[page 58]{brown1994} and since homology commutes with direct limit 
		\cite[Exer.~6, \S5, Chap.~V]{brown1994}, we have $H_n(A,\z)=0$ for $n\geq 2$. Now 
		the claim follows from an easy analysis of the above Lyndon/Hochschild-Serre spectral 
		sequence. The second claim follows from Corollary \ref{2-inv}.
	\end{proof}
	\begin{exa}\label{ZZ-0}
		(i) Since $\z^\times$ act on $\z$ by $(-1).n=-n$, we have $(\z)_{\z^\times}\simeq \z/2$. Moreover, 
		using the structure  of the homology of finite cyclic groups \cite[pages 58-59]{brown1994}, we have 
		$H_k(\z^\times,\z)\simeq\begin{cases}
			0 & \text{if $k$ is odd}\\
			\z/2 & \text{if $k$ is even}
		\end{cases}$.
		Therefore, by the above proposition,
		\begin{align*}
			H_n(\PB_2(\z),\z)\simeq
			\begin{cases}
				H_n(\PT_2(\z),\z) & \text{if $n$ is even}\\
				H_n(\PT_2(\z),\z)\oplus \z/2 & \text{if $n$ is odd}
			\end{cases}
			\simeq 
			\begin{cases}
				\z & \text{if $n=0$}\\
				0 & \text{if $n$ is even.}\\
				\z/2 \oplus \z/2 & \text{if $n$ is odd}
			\end{cases}
		\end{align*}
		
		\par (ii) Let  $A_m:=\z[\frac{1}{m}]$, where $m$ is a square free integer. We calculated $(A_m)_{A_m^\times}$
		in Example \ref{exa-H1}(iv). Now if $2\mid m$, then by Corollary \ref{2-inv}, for any non-negative integer $n$, 
		$H_n(A_m^\times, A_m)=0$. Thus by Proposition \ref{1/2} 
		\[
		H_n(\PB_2(A_m),\z)\simeq H_n(\PT_2(A_m),\z).
		\]
		Let $p$ be an odd prime. Then ${(A_p)}_{A_p^\times}\simeq \z/2$. Note that 
		$A_p^\times\simeq \{\pm 1\} \times \lan p\ran\simeq \z/2 \times \z$. By an easy analysis
		of the spectral sequence
		\[
		{\EE'}_{r,s}^2=H_r(\{\pm 1\}, H_s(\lan p\ran, A_p)) \Longrightarrow H_{r+s}(A_p^\times, A_p)
		\]
		and the calculation of the homology of cyclic groups \cite[Pages 58-59]{brown1994}, one can show 
		that for any $n\geq 0$,
		\[
		H_n(A_p^\times, A_p)\simeq \z/2.
		\]
		Thus by Proposition \ref{1/2}, $H_n(\PB_n(A_p),\z)\simeq H_n(\PT_2(A_p),\z)\oplus \z/2$, $n\geq 1$.
		Finally by the K\"unneth formula applied to
		\[
		H_n(\PT_2(A_p),\z)\simeq H_n(\{\pm 1\} \times \lan p\ran,\z),
		\]
		we obtain
		\[
		H_n(\PB_2(A_p),\z)\simeq 
		\begin{cases}
			\z & \text{if $n=0$} \\
			\z/2 \oplus \z/2 \oplus \z & \text{if $n=1$} \\
			\z/2 \oplus \z/2 & \text{if $n\geq 2$} \\
		\end{cases}.
		\]
	\end{exa}
	
	The following result is known for the inclusion of groups ${\rm T}_2(A)\se {\rm B}_2(A)$.
	
	\begin{prp}\label{B2-T2}
		Let $\psi_q:H_q(\PT_2(A),\z)\arr H_q(\PB_2(A),\z)$ be induced by the natural inclusion $\PT_2(A)\se \PB_2(A)$.
		\par {\rm (i)} Let $A$ be a semilocal ring such that for any maximal ideal $\mmm$, $|A/\mmm|\neq 2,3, 4$. 
		Then $\psi_q$ is isomorphism for $q\leq 2$.
		\par {\rm (ii)} Let $A$ be a semilocal domain such that for any maximal ideal $\mmm$ either $A/\mmm$ is infinite 
		or if $|A/\mmm|=p^d$, then $q<(p-1)d$. Then $\psi_q$ is an isomorphism.
		\par {\rm (iii)} Let $A$ be a semilocal ring such that for any maximal ideal $\mmm$ either $A/\mmm$ is 
		infinite or if $|A/\mmm|=p^d$, then $q<(p-1)d-2$. Then $\psi_q$ is an isomorphism.
	\end{prp}
	\begin{proof}
		This can be proved as in \cite[Section \S2]{mirzaii2017}, which the case of local rings is treated.
	\end{proof}
	
	\section{The second homology of  \texorpdfstring{$\PGE_2$}{Lg}}\label{sec9}
	
	Let $A$ be a commutative ring. Recall that $\WW_A=\{a\in A:a(a-1)\in \aa\}$.
	The differential $\partial_3^Y:Y_3(A) \arr Z_2^{\GE_2}(A^2)\se Y_2(A)$ induces the map
	\[
	\overline{\partial_3^Y}: H_0(\PGE_2(A), Y_3(A^2))\arr \PP(A).
	\]
	We choose $X_a:=({\pmb \infty}, {\pmb 0} ,{\pmb 1}, {\pmb a})$, $a\in \WW_A$, as representatives 
	of the orbits of the generators of $Y_3(A^2)$ and set
	\[
	[a]:=\overline{\partial_3^Y}(X_a)\in \PP(A).
	\]
	
	\begin{prp}\label{d21}
		Let $A$ be a ring such that $Y_\bullet(A^2)$ is exact in dimension $1$. Then under the composite
		\[
		d_{2,1}^2:\PP(A) \overset{d_{2,1}^2}{\larr} H_2(\PB_2(A),\z) \arr H_2(\PT_2(A),\z) \simeq \aa\wedge \aa
		\]
		$[a]\in \PP(A)$ maps to $2\big(a\wedge (1-a)\Big)$.
	\end{prp}
	\begin{proof}
		It is proved in \cite[Lemma 3.2]{mirzaii2011} that ${E'}_{2,1}^2\simeq \PP(A)$ and the composite
		\[
		\PP(A) \overset{{d'}_{2,1}^2}{\larr} \frac{H_2({\rm B}_2(A),\z)}{(\sigma_\ast-\inc_\ast)(H_2({\rm T}_2(A),\z))} 
		\larr \frac{H_2({\rm T}_2(A),\z)}{(\sigma_\ast-\inc_\ast)(H_2({\rm T}_2(A),\z))}
		\]
		is given by 
		\[
		[a]\mapsto {\mtxx{a}{0}{0}{1}}\wedge {\mtxx{1-a}{0}{0}{(1-a)^{-1}}}
		\]
		(see \cite[Lemma 4.1]{mirzaii2011} and its proof). We introduced and studies the spectral sequence 
		\[
		{E'}_{p,q}^1\Rightarrow H_{p+q}(\GE_2(A),\z)
		\]
		in the proof of Lemma \ref{E21}. Observe that in \cite{mirzaii2011} the author works over rings with many 
		units, which satisfy the condition of this proposition. But Lemma 4.1 of \cite{mirzaii2011} (used above) is also valid 
		in our more general setting here. Now from the commutative diagram
		\[
		\begin{tikzcd}
			\PP(A) \ar[r, "{d'}_{2,1}^2"] \ar[d, equal] &\displaystyle \frac{H_2({\rm B}_2(A),\z)}{\im(\sigma_\ast-\inc_\ast)}
			\ar[d, "p_\ast"] \ar[r] & \displaystyle\frac{H_2({\rm T}_2(A),\z)}{\im(\sigma_\ast-\inc_\ast)} \ar[d, "p_\ast"]\\
			\PP(A)  \ar[r, "{d}_{2,1}^2"] & H_2(\PB_2(A),\z) \ar[r] & H_2(\PT_2(A),\z)
		\end{tikzcd}
		\]
		we obtain the desired result.
	\end{proof}
	
	Let $A$ satisfy the condition that $Y_\bullet(A^2)$ is exact in dimension $1$. Then by 
	Lemma~\ref{E21}, $E_{2,1}^2\simeq \PP(A)$.
	We denote the kernel of the differential
	\[
	d_{2,1}^2:\PP(A) \arr H_2(\PB_2(A),\z)
	\]
	with $\BB_E(A)$ and we call it the {\it $\GE_2$-Bloch group} of $A$. Hence we have
	\begin{equation}\label{E21=B}
		E_{2,1}^\infty\simeq \BB_E(A).
	\end{equation}
	
	\begin{cor}\label{H2}
		Let $A$ satisfy the condition that $Y_\bullet(A^2)$ is exact in dimension $1$.
		If $H_k(\PT_2(A),\z)\simeq H_k(\PB_2(A),\z)$ for $k\leq 2$, then 
		we have the exact sequence
		\[
		\frac{\aa \wedge \aa}{\lan 2(a\wedge (1-a)):a\in \WW_A\ran} \arr H_2(\PGE_2(A),\z) \arr \mu_2(A)\arr 1.
		\]
		Moreover, if $Y_\bullet(A^2) \arr\z$ is exact in dimension $<3$, then we have the exact sequence
		\[
		0 \arr \frac{\aa \wedge \aa}{\lan 2(a\wedge (1-a)):a\in \WW_A\ran} \arr H_2(\PGE_2(A),\z) \arr \mu_2(A)\arr 1.
		\]
	\end{cor}
	\begin{proof}
		This follows from Theorem \ref{ES} and Proposition \ref{d21}. The second part follows from the first part and 
		the fact that the natural map 
		\[
		\overline{\partial_3^Y}: H_0(\PGE_2(A), Y_3(A^2))\arr \PP(A)
		\]
		is surjective.
	\end{proof}
	
	\begin{thm}
		Let $A$ be a local domain (local ring) such that $|A/\mmm_A|\neq 2, 3, 4$ 
		($|A/\mmm_A|\neq 2, 3, 4, 5,8,9,16$). Then
		\[
		\begin{array}{c}
			H_2(\PGL_2(A),\z\half)\simeq \K_2(A)\half.
		\end{array}
		\]
	\end{thm}
	\begin{proof}
		Since $|A/\mmm_A|\neq 2, 3, 4$ ($|A/\mmm_A|\neq 2, 3, 4, 5,8,9,16$ for the case of local ring), by 
		Proposition \ref{B2-T2}, 
		\[
		H_k(\PT_2(A),\z)\simeq H_k(\PB_2(A),\z)
		\]
		for $k\leq 2$. Moreover, note that $A$ is a $\GE_2$-ring and by Proposition \ref{<k} the complex 
		$Y_\bullet(A^2)\arr \z$ is exact in dimension $<4$. Thus by the above corollary,
		\[
		\begin{array}{c}
			H_2(\PGL_2(A),\z\half) \simeq \displaystyle\frac{H_2(\aa,\z\half)}{\lan 2(a\wedge (1-a)):a\in \WW_A\ran}.
		\end{array}
		\]
		Now it is easy to see that
		\begin{align*}
			\displaystyle\frac{H_2(\aa,\z\half)}{\lan 2(a\wedge (1-a)):a\in \WW_A\ran} 
			\simeq \frac{S_\z^2(\aa)\half}{\lan a\otimes (1-a):a\in \aa\ran}
			\begin{array}{c}
				\simeq \K_2^M(A)\half \simeq \K_2(A)\half
			\end{array}
		\end{align*}
		(see Theorem \ref{MV}). Recall that $S_\z^2(\aa)\simeq (\aa\otimes_\z \aa)/\lan a\otimes b +b \otimes a:a, b\in \aa\ran$.
	\end{proof}
	
	\begin{exa}\label{ZZ}
		The ring $\z$ is an universal $\GE_2$-ring \cite[Example 6.12]{hutchinson2022}. 
		Since $\z_{\z^\times}\simeq \z/2$ and $H_2(\PB_2(\z),\z)=0$ (Example \ref{ZZ-0}), 
		by Corollary~\ref{H1-uni} we have 
		\[
		H_2(\PGL_2(\z),\z)\simeq \z/2.
		\]
	\end{exa}
	
	\begin{exa}\label{ZZZ}
		Let $p=2,3$. Then $A_p:=\z[\frac{1}{p}]$ is an universal $\GE_2$-ring \cite[Example 6.13]{hutchinson2022}.
		By Example \ref{ZZ-0} we have  
		\[
		H_2(\PB_2(A_p),\z)\simeq 
		\begin{cases}
			\z/2 & \text{if $p=2$}\\
			\z/2\oplus \z/2 & \text{if $p=3$.}
		\end{cases}
		\]
		By Corollary \ref{H1-uni} we have the exact sequence
		\[
		H_3(\PGL_2(A_p),\z) \arr \PP(A_p) \overset{\lambda}{\larr} H_2(\PB_2(A_p),\z) \arr H_2(\PGL_2(A_p),\z) \arr \mu_2(A_p) \arr 1.
		\]
		Since $H_2(\PB_2(A_p),\z)$ is 2-torsion, we have $\lambda([a])=2(a\wedge (1-a))=0$. 
		Coronado and Hutchinson have studied $\PP(A_2)$ in \cite[Corollary 8.42]{C-H2022}. They showed that
		\[
		\begin{array}{c}
			\PP(A_2)\simeq \z/12 \oplus \z/2
		\end{array}
		\]
		and  has generators $[2]$ and $\half$, both of order $12$ and satisfying $2([2] +\half) = 0$. It follows from these results
		that $\PP(A_2) \overset{\lambda}{\larr} H_2(\PB_2(A_2),\z)\simeq \z/2$ is trivial and hence we have the exact sequence
		\[
		0 \arr \z/2 \arr H_2(\PGL_2(A_2),\z)\arr \mu_2(A_2) \arr 1.
		\]
	\end{exa}
	
	\section{The third homology of \texorpdfstring{$\PGE_2$}{Lg} and a Bloch-Wigner type theorem} \label{sec10}
	
	\begin{lem}\label{E12}
		If $A$ satisfies the condition that $Y_\bullet(A^2)$ is exact in dimension $1$, then 
		\[
		E_{1,2}^2\simeq \displaystyle\frac{\aa \wedge \aa}{2(\aa \wedge \aa)}\simeq \GG_A\wedge \GG_A.
		\]
	\end{lem}
	\begin{proof}
		To prove our claim we must show that the image of the differential 
		\[
		d_{2,2}^1: H_2(\PGE_2(A),Z_1^{\GE_2}(A^2))\arr H_2(\PT_2(A),\z)\simeq \aa \wedge \aa
		\]
		is $2(\aa \wedge \aa)$. The morphism of spectral sequences ${E'}_{p,q}^1\arr  E_{p,q}^1$ (see diagram 
		(\ref{morphism}) in the proof of Lemma \ref{E21}) gives us the commutative diagram
		\[
		\begin{tikzcd}
			H_2(\GE_2(A), Z_1^{\GE_2}(A^2)) \ar[r, "{d'}_{2,2}^1"] \ar[d] & H_2({\rm T}_2(A),\z) \ar[r,"{d'}_{1,2}^1"]\ar[d, "p_\ast"] 
			& H_2({\rm B}_2(A),\z)\ar[d]\\
			H_2(\PGE_2(A), Z_1^{\GE_2}(A^2)) \ar[r, "d_{2,2}^1"] & H_2(\PT_2(A),\z) \ar[r, "d_{1,2}^1=0"]       & H_2(\PB_2(A),\z).
		\end{tikzcd}
		\]
		From the  five term exact sequence associated to the Lyndon/Hochschild-Serre spectral of central extension 
		$1 \arr \aa I_2 \arr  \GE_2(A) \arr \PGE_2(A)\arr 1$ \cite[Corollary 6.4, Chp. VII]{brown1994}, with coefficients 
		in $Z_1^{\GE_2}(A^2)$ and Lemma \ref{PE} we obtain the exact sequence
		\[
		H_2(\GE_2(A), Z_1^{\GE_2}(A^2)) \arr H_2(\PGE_2(A), Z_1^{\GE_2}(A^2))\arr H_1(\aa,Z_1^{\GE_2}(A^2))_{\PGE_2(A)} 
		\]
		\[
		\arr H_1(\GE_2(A), Z_1^{\GE_2}(A^2)) \arr \PP(A)\arr 0.
		\]
		Since $\PGE_2(A)$ acts trivially on $\aa I_2$, we have
		\[
		H_1(\aa,Z_1^{\GE_2}(A^2))_{\PGE_2(A)}\simeq H_1(\aa,\z)\otimes_\z Z_1^{\GE_2}(A^2)_{\PGE_2(A)}.
		\]
		In the proof of Lemma \ref{PE}, we have proved that $Z_1^{\GE_2}(A^2)_{\PGE_2(A)}\simeq \z$. Thus 
		\[
		H_1(\aa,Z_1^{\GE_2}(A^2))_{\PGE_2(A)}\simeq \aa.
		\]
		The composite $Y_2(A^2) \overset{\partial_2}{\two} Z_1^{\GE_2}(A^2) \overset{\inc}{\larr} Y_1(A^2)$
		gives us the composite
		\[
		\aa\simeq H_1(\GE_2(A), Y_2(A^2)) \arr H_1(\GE_2(A), Z_1^{\GE_2}(A^2)) \arr H_1(\GE_2(A), Y_1(A^2))\simeq \aa
		\]
		which coincide with the identity map $\id_\aa$. This shows that the natural map
		\[
		\aa \arr  H_1(\GE_2(A), Z_1^{\GE_2}(A^2))
		\]
		appearing in the above five term exact sequence is injective. Therefore the map 
		\[
		H_2(\GE_2(A), Z_1^{\GE_2}(A^2))\arr H_2(\PGE_2(A), Z_1^{\GE_2}(A^2))
		\]
		is surjective, which appears as the left vertical map in the above diagram.
		
		By the K\"unneth formula,
		\[
		H_2({\rm T}_2(A),\z)\simeq H_2(\aa,\z)\oplus H_2(\aa,\z)\oplus \aa\otimes_\z \aa.
		\]
		A direct calculation shows that $\ker({d'}_{1,2}^1)$ is generated by the elements of the
		form 
		\[
		(x,x, a\otimes b-b\otimes a), \ \ \ \ x\in H_2(\aa,\z),\ \ \ a, b\in \aa.
		\]
		It is easy to see that $p_\ast (\ker{d'}_{1,2}^1)=2H_2(\PT(A),\z)$. From the above commutative 
		diagram it follows that $\im(d_{2,2}^1)\se 2H_2(\PT(A),\z)$. Finally it is straightforward to check that 
		\[
		Y:=\!(\bigg[{\mtxx{a}{0}{0}{1}}|{\mtxx{b}{0}{0}{1}}\bigg]\!-\!\bigg[{\mtxx{b}{0}{0}{1}}|{\mtxx{a}{0}{0}{1}}\bigg])
		\!\otimes\!\Big(({\pmb\infty},{\pmb 0})\!+\!({\pmb 0},{\pmb \infty})\Big)\! \in\! H_2(\PGE_2(A),Z_1^{\GE_2}(A^2))
		\]
		and 
		\[
		d_{2,2}^1(Y)=2(a\wedge b).
		\]
		This shows that $\im(d_{2,2}^1)=2H_2(\PT(A),\z)$. The final isomorphism 
		follows from the following lemma applied two $A = \z \arr \z/2 = B$ and $M=\aa$.
	\end{proof}
	
	\begin{lem}[Base change]\label{wedge}
		If $A \arr B$ is a homomorphism of commutative rings and if $M$ is any $A$-module, then the natural map 
		$(\bigwedge_A^n M)\otimes_A B \arr \bigwedge_B^n (M\otimes_A B)$ is an isomorphism.
	\end{lem}
	\begin{proof}
		See \cite[Proposition A2.2]{eisenbud1995}
	\end{proof}
	
	Let $\Aa$ be an abelian group. Let $\sigma_1:\tors(\Aa,\Aa)\arr \tors(\Aa,\Aa)$
	be obtained by interchanging the copies of $\Aa$. This map is induced by
	the involution $\Aa\otimes_\z \Aa \arr \Aa\otimes_\z \Aa$, $a\otimes b\mapsto -b\otimes a$
	\cite[\S2]{mmo2024}. Let $\Sigma_2'=\{1, \sigma'\}$ be the symmetric group of order 2 and 
	consider the following action of this group on $\tors(\Aa,\Aa)$:
	\[
	(\sigma', x)\mapsto -\sigma_1(x).
	\]
	We say that an abelian group $\Aa$ is an
	ind-cyclic group if $\Aa$ is direct limit of its finite cyclic subgroups.
	
	\begin{prp}\label{H3A}
		Let $\Aa$ be an abelian group and $T_\Aa$ be its torsion subgroup. Then 
		\par {\rm (i)} We always have the exact sequence
		\[
		\begin{array}{c}
			0 \arr \bigwedge_\z^3 \Aa\arr H_3(\Aa,\z) \arr \tors(T_\Aa,T_\Aa)^{\Sigma_2'}\arr 0.
		\end{array}
		\]
		\par {\rm (ii)}
		If $T_\Aa$ is an ind-cyclic group, then $\Sigma_2'$ acts trivially on $\tors(T_\Aa,T_\Aa)$  
		and the exact sequence 
		\[
		\begin{array}{c}
			0 \arr \bigwedge_\z^3 \Aa\arr H_3(\Aa,\z) \arr \tors(T_\Aa,T_\Aa) \arr 0
		\end{array}
		\]
		splits naturally.
		\par {\rm (iii)} For any integer $m\in \z$, let $m:\Aa \arr \Aa$ be given by $a\mapsto ma$. Then the map 
		$m_\ast:H_3(\Aa,\z) \arr H_3(\Aa,\z)$ induces multiplication by $m^3$ on $\bigwedge_\z^3 \Aa$ 
		and multiplication by $m^2$ on $\tors(T_\Aa,T_\Aa)^{\Sigma_2'}$.
		
	\end{prp}
	\begin{proof}
		See \cite[Proposition 1.3, Corollary 1.4]{B-E--2024}.
	\end{proof}
	
	\begin{lem}\label{E03}
		If  $H_3(\PT_2(A),\z)\simeq H_3(\PB_2(A),\z)$, then $E_{0,3}^2$  sits in the exact sequence
		\[
		\begin{array}{c}
			\bigwedge_\z^3 \GG_A\arr E_{0,3}^2 \arr \tors(\mu(A),\mu(A))^{\Sigma_2'}\arr 0.
		\end{array}
		\]
	\end{lem}
	\begin{proof}
		By Proposition \ref{H3A}, $H_3(\PT_2(A),\z)$ sits in the exact sequence
		\[
		\begin{array}{c}
			0\arr \bigwedge_\z^3 \aa\arr H_3(\PT_2(A),\z) \arr \tors(\mu(A),\mu(A))^{\Sigma_2'}\arr 0.
		\end{array}
		\]
		Consider the commutative diagram with exact rows
		\[
		\begin{tikzcd}
			0 \ar[r] & \bigwedge_\z^3 \aa \ar[r] \ar[d] & H_3(\PT_2(A),\z) \ar[r]\ar[d, "d_{1,3}^1"] & 
			\tors(\mu(A),\mu(A))^{\Sigma_2'} \ar[r] \ar[d]& 0\\
			0\ar[r]& \bigwedge_\z^3 \aa \ar[r] &  H_3(\PT_2(A),\z)\ar[r] & \tors(\mu(A),\mu(A))^{\Sigma_2'} \ar[r] & 0.
		\end{tikzcd}
		\]
		It is straightforward to see that $d_{1,3}^1$ induces multiplication by $2$ on the left vertical map and 
		$0$ on the right vertical map (use Proposition \ref{H3A} for $m=1$ and $m=-1$). Thus by the Snake lemma 
		we have the exact sequence
		\[
		\begin{array}{c}
			(\bigwedge_\z^3 \aa)/2\arr E_{0,3}^2 \arr \tors(\mu(A),\mu(A))^{\Sigma_2'}\arr 0.
		\end{array}
		\]
		Now the claim follows from this and Lemma \ref{wedge}.
	\end{proof}
	
	Over algebraically closed fields of characteristic zero the following result is called the 
	classical Bloch-Wigner exact sequence.
	
	\begin{prp}[Classical Bloch-Wigner exact sequence]\label{classical}
		Let $F$ be either a quadratically closed field, real closed field or a finite field,
		where $|F|\neq 2,3, 4, 8$. Then we have the exact sequence
		\[
		0\arr \tors(\mu(F),\mu(F)) \arr H_3(\PGL_2(F),\z) \arr \BB_E(F) \arr 0.
		\]
	\end{prp}
	\begin{proof}
		First note that $Y_\bullet(F^2)\arr \z$ is exact in dimension $<4$ (Proposition \ref{<k}). Second,
		by Proposition \ref{B2-T2}, $H_n(\PT_2(F),\z)\simeq H_n(\PB_2(F),\z)$ for $n\leq 3$. Since $|\GG_F|\leq 2$, 
		we have $E_{1,2}^2=0$,  $E_{0,3}^2\simeq \tors(\mu(F),\mu(F))$ (see Lemmas \ref{E12} and \ref{E03}). 
		Now by an easy analysis of the main spectral sequence we obtain the exact sequence
		\[
		\tors(\mu(F),\mu(F)) \arr H_3(\PGL_2(F),\z) \arr \BB_E(F) \arr 0.
		\]
		Let $\of$ be the algebraic closure of $F$. Since $\tors(\mu(F),\mu(F)) \arr \tors(\mu(\of),\mu(\of))$ is injctive,
		it is sufficient to prove the claim for $\of$. The isomorphism $\PSL_2(\of) \overset{\simeq}{\arr} \PGL_2(\of)$
		gives us the morphism of spectral sequences
		\[
		\begin{tikzcd}
			{E''}_{p,q}^1\ar[r,Rightarrow] \ar[d]& H_{p+q}(\PSL_2(\of),\z)\ar[d]\\
			E_{p,q}^1\ar[r,Rightarrow]& H_{p+q}(\PGL_2(\of),\z).
		\end{tikzcd}
		\]
		This morphism gives us the commutative diagram with exact rows
		\[
		\begin{tikzcd}
			0 \ar[r]&\tors(\widetilde{\mu}(\of),(\widetilde{\mu}(\of))\ar[r] \ar[d, "{4.}"] & 
			H_3(\PSL_2(\of),\z)\ar[r]\ar[d, "\simeq"] & \PP(\of)\! \ar[rrr, "{[a]\mapsto \frac{1}{2}(a\wedge (1-a))}"] \ar[d, equal] 
			& & &\! \bigwedge_\z^2\off  \ar[d, "{4.}"]\\
			& \tors(\mu(\of),\mu(\of)) \ar[r] & H_3(\PGL_2(\of),\z) \ar[r] & \PP(\of)\! \ar[rrr, "{[a]\mapsto 2(a\wedge (1-a))}"] 
			& & &\! \bigwedge_\z^2\off
		\end{tikzcd}
		\]
		where $\widetilde{\mu}(\of)=\mu(\of)/\mu_2(\of)$. For the upper exact sequence see \cite[Appendix A]{dupont-sah1982} or
		\cite[Theorem 3.1]{B-E-2024}. Note that the right and the left vertical maps are induced by $\off \overset{(.)^2}{\larr} \off$ 
		and both are isomorphism. Therefore the map 
		\[
		\tors(\mu(\of),\mu(\of)) \arr  H_3(\PGL_2(\of),\z)
		\]
		is injective. This completes the proof of the proposition.
	\end{proof}
	
	\begin{rem}
		Let $F$ be a quadratically closed field. All the groups of the commutative diagram
		\[
		\begin{tikzcd}
			\SL_2(F) \ar[r, two heads] \ar[d, hook] & \PSL_2(F)\ar[d, "\simeq"]\\
			\GL_2(F)           \ar[r, two heads] & \PGL_2(F).
		\end{tikzcd}
		\]
		act on the complex $L_\bullet(F^2)$ (see Section \ref{sec4}).  So from the above diagram we obtain the diagram of morphisms of 
		spectral sequences
		\[
		\begin{tikzcd}
			& {E'''}_{p,q}^1 \ar[dl] \ar[rr, Rightarrow] \arrow[dd] & & H_{p+q}(\SL_2(F),\z) \arrow[dl] \arrow[dd] \\
			{E''}_{p,q}^1 \arrow[rr, Rightarrow, crossing over] \arrow[dd] & & H_{p+q}(\PSL_2(F),\z) \ar[dd] \\
			& {E'}_{p,q}^1 \arrow[dl] \arrow[rr, Rightarrow] & & H_{p+q}(\GL_2(F),\z) \arrow[dl] \\
			E_{p,q}^1 \arrow[rr, Rightarrow] & & H_{p+q}(\PGL_2(F),\z)\arrow[from=uu, crossing over]\\
		\end{tikzcd}
		\]
		By studying the  spectral sequence in the above diagram, we obtain the following commutative diagram with exact rows
		{\tiny
			\[
			\!\!\!\!\!\!\!\!\!\!\!\!\!\!\!\!\!\!\!\!\!\!\!\!\!\!\!\!\!\!\!\!\!\!\!\!\!\!\!\!\!\!\!\!\!\!\!\!\!\!\!\!\!\!\!
			\begin{tikzcd}
				& 0 \ar[r] & \mu(F) \ar[dl, two heads] \arrow[rr] \ar[dd, "{2.}"] & & H_3(\SL_2) \ar[dl, two heads] \ar[dd] \ar[rr] & & 
				\PP(F) \ar[dl, equal] \ar[dd, equal] \ar[rr,"{[a]\mapsto \frac{1}{2}(a\wedge (1-a))}"]& & 
				\bigwedge_\z^2\ff \ar[dl, equal] \ar[dd, "\alpha"] \ar[rr] & & H_2(\SL_2) \arrow[dl, hook] \ar[dd, hook] \ar[r] &  0\\
				0\ar[r] &\widehat{\mu}(F) \ar[rr] \ar["{4.}", dd, two heads] & & H_3(\PSL_2) \ar[rr] \ar[dd, "\simeq"]  & & 
				\PP(F) \ar[dd, equal] \ar[rr,"{\ \ [a]\mapsto \frac{1}{2}(a\wedge (1-a))}"] & & \bigwedge_\z^2\ff \ar[rr] \ar[dd, "4."] \ar[rr] & & 
				H_2(\PSL_2)\ar[dd, "\simeq"] \ar[r]  & \widetilde{\mu}_4(F) \larr 0\\
				%
				& \ \ \ \ \ \ \ \ \ \ 0\ar[r] & \mu(F)\ar[dl, "{2.}"] \ar[rr] & &\widetilde{H}_3(\GL_2) \ar[dl] \arrow[dl] \ar[rr] & &
				\PP(F) \ar[dl, equal] \ar[rr, "{[a]\mapsto (a\wedge (1-a),-a\otimes (1-a))}"] & & \bigwedge_\z^2\ff \oplus S_\z^2(\ff) 
				\ar[dl,"\beta"]\ar[rr]& &  H_2(\GL_2) \ar[dl] \ar[dl] \ar[r]  & 0\\
				0 \ar[r] &\mu(F) \ar[rr] & & H_3(\PGL_2) \ar[rr] & & \PP(F)  \ar[rr,"{\ \ [a]\mapsto 2(a\wedge (1-a))}"] & & 
				\bigwedge_\z^2\ff \ar[rr] && H_2(\PGL_2) \ar[r] & \mu_2(F) \larr 0\\
			\end{tikzcd}
			\]
		}
		where $\alpha(a\wedge b)=2(a\wedge b, a\otimes b)$, 
		$\beta(a\wedge b,c\otimes d)=a\wedge b-c\wedge d$, $\widehat{\mu}(F):=\mu(F)/\mu_4(F)$ and
		$\widetilde{H}_3(\GL_2)$ is the cokernel of the composite
		\[
		\begin{array}{c}
			\bigwedge_\z^3 {\rm T}_2(F) \arr H_3({\rm T}_2(F),\z) \arr H_3(\GL_2(F),\z).
		\end{array}
		\]
		In the above diagram, the exact sequences corresponding to $\SL_2(F)$ and $\PSL_2(F)$ are proved in \cite[App. A]{dupont-sah1982}. 
		For these see also \cite{hutchinson-2013}, \cite{B-E--2023}, \cite{B-E-2024} and \cite{B-E--2024}. For the exact sequence related 
		to $\GL_2(F)$ see \cite{suslin1991} and \cite{mirzaii2011}. The exact sequence related to $\PGL_2(F)$ is the topic of the current 
		article (but also see \cite[App. C, (C.3)]{parry-sah1983}). The maps on the right vertical square sit in the following exact sequences
		\[
		0 \arr H_2(\SL_2(F),\z) \arr H_2(\PSL_2(F),\z)\arr \mu_2(F) \arr 1,
		\]
		\[
		0 \arr H_2(\SL_2(F),\z) \arr H_2(\GL_2(F),\z) \overset{\det}{\larr} H_2(\ff,\z) \arr 0,
		\]
		\[
		\ff\otimes_\z \ff \arr H_2(\GL_2(F),\z) \arr H_2(\PGL_2(F),\z) \arr \ff \overset{(.)^2}{\larr} \ff \arr \GG_F \arr 1.
		\] 
		These exact sequences can be obtained by analysis of the Lyndon/Hochschild-Serre Spectral sequences associated to the group 
		extensions $1 \arr \mu_2(F) \arr \SL_2(F) \arr \PSL_2(F) \arr 1$, $1 \arr \SL_2(F) \arr \GL_2(F) \overset{\det}{\arr} \ff \arr 1$
		and $1 \arr \ff I_2 \arr \GL_2(F) \arr \PGL_2(F) \arr 1$, respectively.
	\end{rem}
	
	\begin{thm}\label{main-thm}
		Let $A$ be a domain satisfying the condition that $Y_\bullet(A)$ is exact in dimension $1$ and
		$H_3(\PT_2(A),\z)\simeq H_3(\PB_2(A),\z)$. Then we have the sequence
		\[
		0 \arr \tors(\mu(A),\mu(A)) \arr H_3(\PGE_2(A),\z) \arr \BB_E(A) \arr 0,
		\]
		which is exact at every term except possibly at the term $H_3(\PGL_2(A), \z)$, where the
		homology of the complex is annihilated by $4$.
	\end{thm}
	\begin{proof}
		The main spectral sequence gives us a filtration
		\[
		0\se F_0H_3\se F_1H_3\se F_2H_3 \se F_3H_3=H_3(\PGE_2(A),\z),
		\]
		where $E_{p,3-p}^\infty\simeq F_pH_3/F_{p-1}H_3$, $0\leq p\leq 3$. Since $E_{3,0}^1=0$, we have $F_2H_3=F_3H_3$.
		From $E_{2,1}^\infty\simeq \BB_E(A)$ (see (\ref{E21=B})) we obtain the exact sequence
		\[
		0\arr F_1H_3 \arr H_3(\PGE_2(A),\z) \arr \BB_E(A) \arr 0.
		\]
		By Lemma \ref{E03} and Proposition \ref{H3A}(ii), $E_{0,3}^2$  sits in the exact sequence
		\[
		\begin{array}{c}
			\bigwedge_\z^3 \GG_A \arr E_{0,3}^2 \arr \tors(\mu(A),\mu(A))\arr 0.
		\end{array}
		\]
		Let $\TT_2=\im(\bigwedge_\z^3 \GG_A \arr E_{0,3}^2)$.
		From the natural inclusion $\mu(A)\harr \aa\simeq \PT_2(A)$, we obtain the diagram with exact rows
		\[
		\begin{tikzcd}
			&  & H_3(\mu(A),\z) \ar[r, "\simeq"]\ar[d] & \tors(\mu(A),\mu(A)) \ar[d, equal]&\\
			0 \ar[r]& \TT_2 \ar[r] &  E_{0,3}^2 \ar[r] & \tors(\mu(A),\mu(A)) \ar[r] & 0
		\end{tikzcd}
		\]
		This shows that the bottom exact sequence splits and thus we have the exact sequence 
		\[
		0 \arr \tors(\mu(A),\mu(A)) \arr E_{0,3}^2 \arr \TT_2 \arr 0.
		\]
		Now from the surjective map $E_{0,3}^2 \two E_{0,3}^\infty\simeq F_0H_3$, 
		we obtain an exact sequence of the form
		\[
		\tors(\mu(A),\mu(A)) \arr F_0H_3 \arr \TT_2' \arr 0,
		\]
		where $\TT_2'$ is a $2$-torsion group. Let $\alpha$ be the composite $\tors(\mu(A),\mu(A)) \arr F_0H_3\harr F_1H_3$.
		From the commutative diagram
		\[
		\begin{tikzcd}
			& \tors(\mu(A),\mu(A)) \ar[r, equal] \ar[d] & \tors(\mu(A),\mu(A)) \ar[d, "\alpha"] & &\\
			0 \ar[r] & F_0H_3\ar[r] &  F_1H_3 \ar[r] & E_{1,2}^\infty \ar[r] & 0
		\end{tikzcd}
		\]
		we obtain the exact sequence 
		\[
		0 \arr \TT_2'\arr \coker(\alpha) \arr E_{1,2}^\infty \arr 0. 
		\]
		But by Lemma \ref{E12}, $E_{1,2}^\infty\simeq \GG_A\wedge \GG_A$. This implies that $\TT_4:=\coker(\alpha)$
		is a $4$-torsion group. To complete the proof of the theorem we need to prove that the composite
		\[
		\tors(\mu(A),\mu(A)) \overset{\alpha}{\larr} F_1H_3 \harr  H_3(\PGE_2(A),\z)
		\]
		is injective. Let $F$ be the quotient field of $A$ and $\overline{F}$ the algebraic closure of $F$.
		By Proposition \ref{classical}, we have the classical Bloch-Wigner exact sequence
		\[
		0\arr \tors(\mu(\overline{F}),\mu(\overline{F}))\arr H_3(\PGL_2(\overline{F}),\z) \arr \BB_E(\overline{F})\arr 0.
		\]
		Now from the commutative diagram 
		\[
		\begin{tikzcd}
			& \tors(\mu(A),\mu(A)) \ar[r] \ar[d, hook] & H_3(\PGE_2(A),\z)\ar[d] & & \\
			0\ar[r]& \tors(\mu(\overline{F}),\mu(\overline{F})) \ar[r] &  H_3(\PGL_2(\overline{F}),\z) \ar[r] & 
			\BB_E(\overline{F}) \ar[r] & 0
		\end{tikzcd}
		\]
		we obtain the injectivity of the map $\tors(\mu(A),\mu(A)) \arr H_3(\PGE_2(A),\z)$. This completes the 
		proof of the theorem.
	\end{proof}
	
	\begin{cor}
		Let $A$ be a semilocal domain such that for any maximal ideal $\mmm$,  $|A/\mmm|\neq 2, 3, 4, 8$. 
		Then we have the Bloch-Wigner exact sequence
		\[
		\begin{array}{c}
			0\arr \tors(\mu(A),\mu(A))\half \arr H_3(\PGL_2(A),\z\half) \arr \BB(A)\half \arr 0.
		\end{array}
		\]
		Moreover, $H_3(\PGL_2(A),\z\half)\simeq \K_3^\ind(A)\half$.
	\end{cor}
	\begin{proof}
		First note that $A$ is a $\GE_2$-ring and thus $\GE_2(A)=\GL_2(A)$. 
		Second, by Proposition~\ref{B2-T2}, $H_n(\PT_2(A),\z)\simeq H_n(\PB_2(A),\z)$ for $n\leq 3$. 
		Consider the commutative diagram with exact rows
		\[
		\begin{tikzcd}
			\PP(A) \ar[r, "\lambda"] \ar[d, equal] & S_\z^2(\aa)\ar[d, "\gamma"]\\
			\PP(A)           \ar[r, "{d}_{2,1}^2"] & H_2(\PT_2(A),\z),
		\end{tikzcd}
		\]
		where 
		$\gamma(\overline{a\otimes b})=2{\mtxx{a}{0}{0}{1}}\wedge {\mtxx{b}{0}{0}{1}}$.
		Note that $\BB(A):=\ker(\lambda)$ is the Bloch group of $A$ (see Section \ref{sec1}).
		Since $S_\z^2(\aa)\half\simeq H_2(\PT_2(A),\z\half)$, we have $\BB(A)\half\simeq \BB_E(A)\half$. Now the first 
		claim follows from the above theorem. 
		
		The natural map $\GL_2(A) \harr \PGL_2(A)$, induces  the commutative diagram with exact rows
		\[
		\begin{tikzcd}
			0 \ar[r] & \tors(\mu(A),\mu(A)\half \ar[r] \ar[d, "\simeq"] & H_3(\SL_2(A),\z\half)_\aa \ar[r]\ar[d] & 
			\BB(A)\half \ar[r] \ar[d, "\simeq"]& 0\\
			0\ar[r]& \tors(\mu(A),\mu(A)\half \ar[r] &  H_3(\PGL_2(A),\z\half)\ar[r] & \BB_E(A)\half \ar[r] & 0
		\end{tikzcd}
		\]
		Note that the first exact sequence can be proved as \cite[Corollary 5.4]{mirzaii2011}.
		By the Snake lemma, $H_3(\SL_2(A),\z\half)_\aa\simeq  H_3(\PGL_2(A),\z\half)$. Now the second claim follows 
		from this and the isomorphism $H_3(\SL_2(A),\z\half)_\aa\simeq  \K_3^\ind(A)\half$ (see Theorem \ref{kind-sl2}).
	\end{proof}
	
	\begin{prp}\label{localfield}
		For any non-dyadic local field $F$ we have the exact sequence
		\[
		0\arr \tors(\mu(F),\mu(F))^\thickapprox \arr H_3(\PGL_2(F),\z) \arr \BB_E(F) \arr 0,
		\]
		where $\tors(\mu(F),\mu(F))^\thickapprox$ is an extension of $\z/2$ by $\tors(\mu(F),\mu(F))$.
	\end{prp}
	\begin{proof}
		First note that $\GG_F\simeq \z/2 \times \z/2$ (see \cite[Theorem 2.2, Chap. VI]{lam2005}).
		Thus $E_{1,2}^\infty\simeq \GG_F\wedge \GG_F\simeq \z/2$. Moreover, since $\bigwedge_\z^3 \GG_F=0$, we have
		$E_{0,3}^2\simeq \tors(\mu(F),\mu(F))$. Now by an easy analysis of the main spectral sequences
		we obtain the exact sequences
		\[
		0\arr \KK \arr H_3(\PGL_2(F),\z) \arr \BB_E(F) \arr 0,
		\]
		\[
		\tors(\mu(F),\mu(F)) \arr \KK \arr \z/2  \arr 0.
		\]
		The injectivity of $\tors(\mu(F),\mu(F)) \arr \KK$ follows from Theorem \ref{main-thm}.
		This completes the proof of the proposition.
	\end{proof}
	
	\begin{prp}
		{\rm (i)} For the ring of integers $\z$, $\BB_E(\z)=\PP(\z)$ and  we have the exact sequence
		\[
		0 \arr \tors(\mu(\z),\mu(\z))\oplus \z/2 \arr H_3(\PGL_2(\z),\z) \arr \BB_E(\z)\arr 0.
		\]
		{\rm (ii)} For the ring $\z\half$, $\BB_E(\z\half)=\PP(\z\half)$ and we have the exact sequence
		\[
		\begin{array}{c}
			0 \arr \tors(\mu(\z\half),\mu(\z\half))^\thickapprox \arr H_3(\PGL_2(\z\half),\z) \arr \BB_E(\z\half)\arr 0,
		\end{array}
		\]
		where $\tors(\mu(\z\half),\mu(\z\half))^\thickapprox$ is an extension of $\z/2$ by $\tors(\mu(\z\half),\mu(\z\half))$.
	\end{prp}
	\begin{proof}
		(i) First observe that  $E_{0,2}^1=H_2(\PB_2(\z),\z)\simeq H_2(\PT_2(\z),\z)=0$ (see Example \ref{ZZ-0}). 
		Thus $\BB_E(\z)=\PP(\z)$. Coronado and Hutchinson have proved that 
		$\PP(\z)\simeq \z/6$ \cite[Corollary 8.20]{C-H2022} and is generated by
		\[
		C_\z:= ({\pmb 1}, {\pmb 0}, {\pmb \infty}) + ({\pmb 1}, {\pmb \infty}, {\pmb 0})
		- ({\pmb 0}, {\pmb \infty}, {\pmb -1})-({\pmb \infty}, {\pmb 0}, {\pmb -1}).
		\]
		Moreover, by Example \ref{ZZ-0} and Proposition \ref{H3A}
		\[
		H_3(\PB_2(\z),\z)\simeq H_3(\PT_2(\z),\z) \oplus \z/2\simeq \tors(\mu(\z),\mu(\z))\oplus \z/2. 
		\]
		Note that $E_{1,2}^1=H_2(\PT_2(\z), \z)=0$. Now by an easy analysis of the main spectral sequence 
		we obtain the exact sequence
		\begin{equation}\label{BW-PGL2(Z)}
			\tors(\mu(\z),\mu(\z))\oplus \z/2 \arr H_3(\PGL_2(\z),\z) \arr \PP(\z)\arr 0.
		\end{equation}
		From this it follows that
		\[
		|H_3(\PGL_2(\z),\z)|\leq 24.
		\]
		On the other hand, we know that $\PGL_2(\z)$ is the free product with amalgamation of the dihedral group $D_2$ 
		of order 4 and the dihedral group $D_3$ of order 6 amalgamated along the subgroup $D_1$ of order $2$:
		\begin{equation}\label{PGL2Z}
			\PGL_2(\z)\simeq D_2\ast_{D_1} D_3,
		\end{equation}
		(see the proof of \cite[Lemma 2]{RT2021}).
		Note that $D_1\simeq \z/2$, $D_2\simeq \z/2\times \z/2$ and $D_3\simeq S_3$. It is easy to see that
		\[
		H_1(D_1,\z)\simeq \z/2, \ \ \ H_2(D_1,\z)=0, \ \ \ H_3(D_1,\z)\simeq \z/2,
		\]
		\[
		H_1(D_2,\z)\simeq (\z/2)^2, \ \ \ H_2(D_2,\z)\simeq \z/2, \ \ \ H_3(D_2,\z)\simeq (\z/2)^3,
		\]
		\[
		H_1(D_3,\z)\simeq \z/2, \ \ \ H_2(D_3,\z)\simeq \z/2, \ \ \ H_3(D_3,\z)\simeq \z/6.
		\] 
		From the Mayer-Vietoris exact sequence associated to (\ref{PGL2Z})  \cite[\S9, Chap. VII]{brown1994} and 
		the above information about the homology groups of $D_1$, $D_2$ and $D_3$ we obtain the exact sequence
		\[
		\z/2 \arr (\z/2)^3 \oplus \z/6 \arr H_3(\PGL_2(\z),\z) \arr 0.
		\]
		Clearly from this it follows that $|H_3(\PGL_2(\z),\z)|\geq 24$. Therefore $H_3(\PGL_2(\z),\z)$
		has 24 elements. Thus the left hand side map of the exact sequence (\ref{BW-PGL2(Z)}) must be injective.
		This proves our claim.
		
		(ii) Let $A_2=\z\half$. Then by Lemma \ref{1/2},  for any $n\geq 0$ we have
		\[
		H_n(\PB_2(A_2),\z)\simeq H_n(\PT_2(A_2),\z).
		\]
		Since $\GG_{A_2}\simeq \z/2\times\z/2$, we have $E_{1,2}^\infty\simeq \z/2$ (Lemma \ref{E12}) and 
		$E_{0,3}^\infty\simeq \tors(\mu(A_2),\mu(A_2))$ (Lemma \ref{E03}). Now from the main spectral sequence 
		we obtain the exact sequences
		\[
		0 \arr \KK\arr H_3(\PGL_2(A_2),\z) \arr \BB_E(A_2)\arr 0,
		\]
		\[
		\tors(\mu(A_2),\mu(A_2)) \arr \KK \arr \z/2\arr 0.
		\]
		As in the proof of Theorem \ref{main-thm}, one can show that the map
		\[
		\tors(\mu(A_2),\mu(A_2)) \arr \KK \se H_3(\PGL_2(A_2),\z)
		\]
		is injective. This completes the proof of the proposition.
	\end{proof}
	
	\begin{rem}
		Let $A$ be a domain. Then, up to isomorphism, there are at most two extension of $\z/2$ by $\tors(\mu(A),\mu(A))$:
		the split and the non-split extensions. This follows from the isomorphism
		\[
		\exts(\z/2, \tors(\mu(A),\mu(A)))\simeq\begin{cases}
			0 & \text{if $\mu_{2^\infty}(F)$ is infinite or $\char(F)=2$} \\
			\z/2 & \text{if $\mu_{2^\infty}(F)$ is finite and $\char(F)\neq2$,}
		\end{cases}
		\]
		where $F$ is the quotient field of $A$ and $\mu_{2^\infty}(F)$ is the 2-power roots of unity in $F$. We believe the extension 
		\[
		0 \arr \tors(\mu(A),\mu(A)) \arr \tors(\mu(A),\mu(A))^\thickapprox \arr \z/2 \arr 0
		\]
		appearing in the above two propositions is the non-split extension. But at the moment we do not know how to prove this.
	\end{rem}

\end{document}